\documentclass[a4paper,10pt]{article}

\usepackage[a4paper,left=3cm,right=3cm,top=2.5cm,bottom=2.5cm]{geometry}
\usepackage{hyperref}
\usepackage{subfig}
\usepackage{pgfplots}
\usepackage{pgfplotstable}
\usepackage{mathtools}
\usepackage{multicol}
\usepackage{comment}
\usepackage{booktabs}
\pgfplotsset{compat=1.5}
\usepackage{amssymb}
\usepackage{url}
\usepackage{bm}
\usepackage{tikz} 
\usetikzlibrary{shapes,arrows,positioning,chains,fit,arrows.meta}

\usepackage{pdfsync}
\usepackage{float}
\usepackage{tabularx}
\usepackage{enumerate}
\usepackage{array}
\usepackage{xspace}
\usepackage{siunitx}
\usepackage{authblk}

\usepackage{amsthm}
\usepackage{amsmath}
\usepackage{stmaryrd}
\usepackage{algorithm}
\usepackage[noend]{algpseudocode}
\usepackage{graphicx}
\usepackage{booktabs}
\usepackage{cleveref}

\DeclareMathOperator*{\argmax}{arg\,max}

\newcommand{\fig}{\text{Fig.~}}
\newcommand{\tab}{\text{Tab.~}}
\newcommand{\eq}{\text{Eq.~}}
\newcommand{\pb}{\text{Problem~}}
\newcommand{\sez}{\text{Sec.~}}
\newcommand{\alg}{\text{Algorithm~}}

\begin{document}

\title{A digital twin framework for civil engineering structures}

\author[1]{Matteo~Torzoni\footnote{matteo.torzoni@polimi.it}}
\author[2]{Marco~Tezzele\footnote{marco.tezzele@austin.utexas.edu}}
\author[1]{Stefano~Mariani\footnote{stefano.mariani@polimi.it}}
\author[3]{Andrea~Manzoni\footnote{andrea1.manzoni@polimi.it}}
\author[2]{Karen~E.~Willcox\footnote{kwillcox@oden.utexas.edu}}

\affil[1]{Dipartimento di Ingegneria Civile e Ambientale, Politecnico di Milano, Piazza L. da Vinci 32, 20133 Milan, Italy}
\affil[2]{Oden Institute for Computational Engineering and Sciences, University of Texas at Austin, Austin, 78712, TX, United States}
\affil[3]{MOX, Dipartimento di Matematica, Politecnico di Milano, Piazza L. da Vinci 32, 20133 Milan, Italy}

\maketitle

\begin{abstract}
The digital twin concept represents an appealing opportunity to advance condition-based and predictive maintenance paradigms for civil engineering systems, thus allowing reduced lifecycle costs, increased system safety, and increased system availability. This work proposes a predictive digital twin approach to the health monitoring, maintenance, and management planning of civil engineering structures. The asset-twin coupled dynamical system is encoded employing a probabilistic graphical model, which allows all relevant sources of uncertainty to be taken into account. In particular, the time-repeating observations-to-decisions flow is modeled using a dynamic Bayesian network. Real-time structural health diagnostics are provided by assimilating sensed data with deep learning models. The digital twin state is continually updated in a sequential Bayesian inference fashion. This is then exploited to inform the optimal planning of maintenance and management actions within a dynamic decision-making framework. A preliminary offline phase involves the population of training datasets through a reduced-order numerical model and the computation of a health-dependent control policy. The strategy is assessed on two synthetic case studies, involving a cantilever beam and a railway bridge, demonstrating the dynamic decision-making capabilities of health-aware digital twins.
\end{abstract}

\tableofcontents

\section{Introduction}

The optimal management of deteriorating structural systems is an important challenge in modern engineering. In particular, the failure or non-optimized maintenance planning of civil structures may entail high safety, economic, and social costs. Within this context, enabling a digital twin (DT) perspective for structural systems that are critical for either safety or operative reasons, is crucial to allow for condition-based or predictive maintenance practices, in place of customarily employed time-based ones. Indeed, having an up-to-date digital replica of the physical asset of interest can yield several benefits spanning its entire lifecycle, including performance and health monitoring, as well as maintenance, inspection, and management planning~\cite{art:DT_review}.

 The DT concept~\cite{tuegel2011reengineering, tuegel2012airframe, grieves2017digital, rasheed2020digital, art:review_1} has been recently applied to several fields for operational monitoring, control, and decision support, including structural health monitoring (SHM) and predictive maintenance~\cite{art:review_2,art:review_3}, additive manufacturing~\cite{art:phua2022digital}, smart cities~\cite{art:review_6}, urban sustainability~\cite{art:tzachor2022potential}, and railway systems management~\cite{arcieri2023bridging}. It allows for a personalized characterization of a physical asset, in the form of computational models and parameters of interest, that evolves over time and is kept synchronized with its physical counterpart by means of data-collecting devices. Within a civil SHM framework, such a twinning perspective can be enabled by the assimilation of data through data-driven structural health diagnostics (from physical to digital), possibly accommodating the quantification and propagation of relevant uncertainties related to, e.g., measurement noise, modeling assumptions, environmental and operational variabilities~\cite{art:kapteyn2020data,art:DT_Cueto,art:Sisson_2, art:moya2020physically, art:cherifi2022simulations}. The resulting updated digital state should then enable prediction of the physical system evolution, as well as inform optimal planning of maintenance and management actions (from digital to physical).

\begin{figure}[h!]
\begin{centering}
\includegraphics[width=1\textwidth]{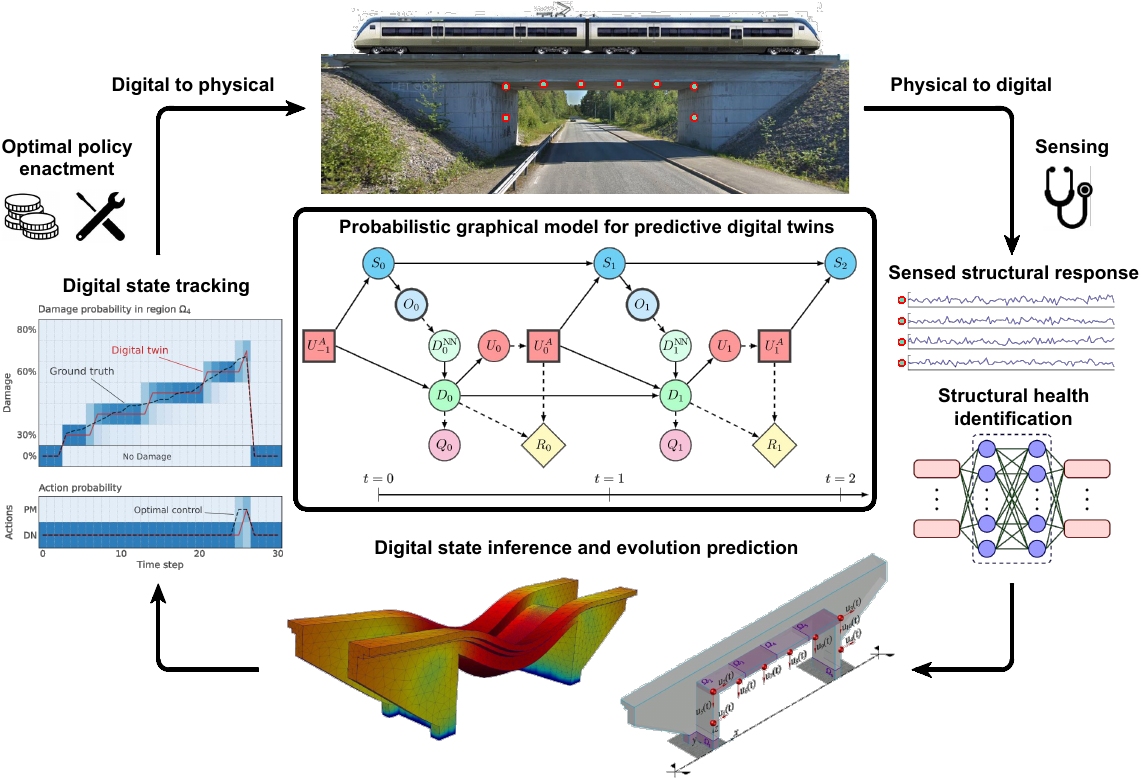}
\caption{Predictive digital twin framework for civil engineering structures: graphical abstraction of the end-to-end information flow enabled by the probabilistic graphical model.\label{fig:graph_abs}}
\end{centering}
\end{figure}

In this work, we propose a DT framework for civil engineering structures. The overall computational strategy is based upon a probabilistic graphical model (PGM) inspired by the foundational model proposed in~\cite{art:kapteyn2021probabilistic}, which provides a general framework to carry out data assimilation, state estimation, prediction, planning, and learning. Formally, such a PGM is a dynamic Bayesian network with the addition of decision nodes, i.e., a dynamic decision network~\cite{book:2009probabilistic,book:2010artificial}. This is employed to encode the end-to-end
information flow, from physical to digital through assimilation and inference, and back to the physical asset in the form of informed control actions. A graphical abstraction of the proposed DT strategy is depicted in \fig\ref{fig:graph_abs}. The figure shows a \textit{physical-to-digital} information flow and a \textit{digital-to-physical} information flow.
These bi-directional information flows repeat indefinitely over time. In particular, we have:

\begin{itemize}
    \item \textit{From physical to digital.} Structural response data are gathered from the physical system and assimilated with deep learning (DL) models, see e.g.,~\cite{art:Torzoni_DML,art:Giglioni}, to estimate the current structural health in terms of presence, location, and severity of structural damage. To solve this inverse problem, we refer to vibration-based SHM techniques, see e.g.,~\cite{art:Torzoni_MF,proc:Torzoni_EWSHM,proc:Rosafalco2022,art:Azam_Mariani}, which exploit the aforementioned collected data, such as displacement or acceleration time histories. This first estimate of the digital state is then employed to estimate an updated digital state, according to control-dependent transition dynamics models describing how the structural health is expected to evolve.
    \item \textit{From digital to physical.} The updated digital state is exploited to predict the future evolution of the physical system and the associated uncertainty, thereby enabling predictive decision-making about maintenance and management actions feeding back to the physical system.
    \item \textit{Offline learning phase.} The DT setup considered in this work takes advantage of a preliminary offline learning phase. This phase involves the training of the DL models underlying the structural health identification, and learning the control policy to be applied at each time step of the online phase. The DL models are trained in a supervised fashion, with labeled data pertaining to specific damage conditions generated by exploiting physics-based numerical models. To efficiently assemble a training dataset representative of potential damage and operational conditions the structure might undergo during its lifetime, we exploit a reduced-order modeling strategy for parametrized systems relying on the reduced basis method~\cite{book:RB}. The health-dependent control policy is also computed offline, by maximizing the expected future rewards for the planning problem induced by the PGM. 
\end{itemize}

The elements of novelty that characterize this work are the following: (i) the adaptation of the PGM digital twinning framework to the health monitoring, maintenance, and management planning of civil engineering structures;
(ii) the assimilation of vibration response data is carried out by exploiting DL models, which allow automated selection and extraction of optimized damage-sensitive features and real-time assessment of the structural state. This work shows how to incorporate in the DT framework high-dimensional multivariate time series describing the sensor measurements, while tracking the associated uncertainties. The proposed computational framework is made available in the public repository \texttt{digital-twin-SHM}~\cite{Repo}. The code implements the PGM framework as a dynamic decision network. This enables us to easily specify the graph topology from a few time slices, and then unroll it for any number of time steps in the future. 

The remainder of this paper is organized as follows. In \sez\ref{sec:methodology}, we describe the proposed DT framework. In \sez\ref{sec:experiments}, the computational procedure is assessed on two test cases, respectively related to a cantilever beam and a railway bridge. Conclusions and future developments are drawn in \sez\ref{sec:conclusions}.

\section{Predictive digital twins using physics-based models and machine learning}
\label{sec:methodology}

In this section, we describe the methodology characterizing our DT framework in terms of the PGM encoding the asset-twin coupled dynamical system in \sez\ref{sec:pgm}; the population of training datasets exploiting physics-based numerical models in \sez\ref{sec:numerical_models}; and the DL models underlying the structural health identification in \sez\ref{sec:DL_models}.

\subsection{Probabilistic graphical model for predictive digital twins}
\label{sec:pgm}

The digital twin assimilates vibration recordings shaped as multivariate time series $\mathbf{U}(\boldsymbol{\mu})=[\mathbf{u}_1(\boldsymbol{\mu}),\ldots,\mathbf{u}_{N_u}(\boldsymbol{\mu})]\in\mathbb{R}^{L\times N_u}$, consisting of $N_u$ time series made of $L$ sensor measurements equally spaced in time, for instance in terms of accelerations or displacements. The vector $\boldsymbol{\mu} \in\mathbb{R}^{N_\text{par}}$ comprises the parameters representing the operational, damage, and (possibly) environmental conditions. Each recording refers to a time interval $(0,T)$, within which measurements are recorded with a sampling rate $f_\text{s}$. For the problem settings we consider, the interval $(0,T)$ is short enough for the operational, environmental, and damage conditions to be considered time-invariant, yet long enough not to compromise the identification of the structural behavior.

The PGM that defines the elements comprising the asset-twin coupled dynamical system, and mathematically describes the relevant interactions via observed data and control inputs, is the dynamic decision network sketched in \fig\ref{fig:graph}. Circle nodes in the graph denote random variables at discrete times, square nodes denote actions, and diamond nodes denote the objective function. Bold outlines denote observed quantities, while thin outlines denote estimated quantities. The directed acyclic structure of the PGM encodes the assumed conditional dependencies. Edges in the graph represent dependencies between random variables. Solid edges represent the variables' dependencies encoded via conditional probability distributions, while dashed edges represent the variables' dependencies encoded via deterministic functions. 

We consider a non-dimensional time discretization, and denote discrete time steps by $t$. The physical time duration between successive time steps may vary depending on the application, and is governed by the update frequency of the DT via data assimilation, so that the DT is updated once per time step. Thanks to the modeled conditional dependencies between random variables, the graph topology is specified from the first two time steps, and can then be unrolled for any number of time steps.

 \tikzstyle{P_node} = [circle,draw=darkgray,thick,align=center,minimum size=0.8cm,fill=cyan!50]
   \tikzstyle{O_node} = [circle,draw=darkgray,line width=.6mm,align=center,minimum size=0.8cm,fill=cyan!20]
 \tikzstyle{D_node} = [circle,draw=darkgray,thick,align=center,minimum size=0.8cm,fill=green!70!cyan!30]
  \tikzstyle{D_NN_node} = [circle,draw=darkgray,thick,align=center,minimum size=0.8cm,fill=green!60!cyan!15]
  \tikzstyle{U_node_prob} = [circle, draw=darkgray, thick, minimum size=0.8cm,fill=red!40]
  \tikzstyle{U_node_act} = [rectangle, draw=darkgray, line width=.6mm, minimum width=0.8cm, minimum height=0.8cm,fill=red!40]
  \tikzstyle{Q_node} = [circle,draw=darkgray,thick,align=center,minimum size=0.8cm,fill=magenta!30]
  \tikzstyle{R_node} = [rectangle,draw=darkgray,thick,align=center,minimum width=0.8cm, minimum height=0.8cm,fill=yellow!30,rotate=45]
  
\begin{figure}[h!]
\center
\begin{tikzpicture}[scale=.85, every node/.style={scale=1.}]

\node [U_node_act] (U_A_-1) at (-3.8,1.5) {};
\node [] at (-3.8,1.5) {$U^A_{-1}$};
\node [P_node] (P_0) at (-2,4) {};
\node [] () at (-2,4) {$S_{0}$};
\node [D_NN_node] (D_NN_0) at (0,1.5) {};
\node [] at (0,1.5) {$D^{\text{N\hspace{-1px}N}}_{0}$};
\node [O_node] (O_0) at (-1,2.75) {};
\node [] at (-1,2.75) {$O_{0}$};
\node [D_node] (D_0) at (0,0) {};
\node [] at (0,0) {$D_{0}$};
\node [U_node_prob] (U_P_0) at (1.5,1.5) {};
\node [] at (1.5,1.5) {$U_{0}$};
\node [U_node_act] (U_A_0) at (3,1.5) {};
\node [] at (3,1.5) {$U^A_{0}$};
\node [Q_node] (Q_0) at (0,-1.5) {};
\node [] at (0,-1.5) {$Q_{0}$};
\node [R_node] (R_0) at (3,-1.5) {};
\node [] at (3,-1.5) {$R_{0}$};

\node [P_node] (P_1) at (5,4) {};
\node [] () at (5,4) {$S_{1}$};
\node [D_NN_node] (D_NN_1) at (7,1.5) {};
\node [] at (7,1.5) {$D^{\text{N\hspace{-1px}N}}_{1}$};
\node [D_node] (D_1) at (7,0) {};
\node [] at (7,0) {$D_{1}$};
\node [O_node] (O_1) at (6,2.75) {};
\node [] at (6,2.75) {$O_{1}$};
\node [U_node_prob] (U_P_1) at (8.5,1.5) {};
\node [] at (8.5,1.5) {$U_{1}$};
\node [U_node_act] (U_A_1) at (10,1.5) {};
\node [] at (10,1.5) {$U^A_{1}$};
\node [Q_node] (Q_1) at (7,-1.5) {};
\node [] at (7,-1.5) {$Q_{1}$};
\node [R_node] (R_1) at (10,-1.5) {};
\node [] at (10,-1.5) {$R_{1}$};

\node [P_node] (P_2) at (12,4) {};
\node [] () at (12,4) {$S_{2}$};

\node [] (t0) at (-2.14,-3) {};
\node [] () at (-2,-3) {$|$};
\node [] () at (-2,-2.5) {$t=0$};
\node [] () at (5,-3) {$|$};
\node [] () at (5,-2.5) {$t=1$};
\node [] () at (12,-3) {$|$};
\node [] () at (12,-2.5) {$t=2$};
\node [] (t2) at (13,-3) {};

\draw[-latex,thick,black] (P_0) to (P_1);
\draw[-latex,thick,black] (P_1) to (P_2);

\draw[-latex,thick,black] (U_A_-1.north east) to (P_0);
\draw[-latex,thick,black] (P_0) to (O_0);
\draw[-latex,dashed,thick,black] (O_0) to (D_NN_0);
\draw[-latex,thick,black] (U_A_-1.east) to (D_0);
\draw[-latex,thick,black] (D_NN_0) to (D_0);
\draw[-latex,thick,black] (D_0) to (U_P_0);
\draw[-latex,dashed,thick,black] (U_P_0) to (U_A_0);
\draw[-latex,dashed,thick,black] (U_A_0) to (R_0);
\draw[-latex,dashed,thick,black] (D_0) to (R_0);
\draw[-latex,dashed,thick,black] (D_0) to (Q_0);

\draw[-latex,thick,black] (D_0) to (D_1);
\draw[-latex,thick,black] (U_A_0.north east) to (P_1);
\draw[-latex,thick,black] (P_1) to (O_1);
\draw[-latex,dashed,thick,black] (O_1) to (D_NN_1);
\draw[-latex,thick,black] (U_A_0.east) to (D_1);
\draw[-latex,thick,black] (D_NN_1) to (D_1);
\draw[-latex,thick,black] (D_1) to (U_P_1);
\draw[-latex,dashed,thick,black] (U_P_1) to (U_A_1);
\draw[-latex,dashed,thick,black] (U_A_1) to (R_1);
\draw[-latex,dashed,thick,black] (D_1) to (R_1);
\draw[-latex,dashed,thick,black] (D_1) to (Q_1);

\draw[-latex,thick,black] (U_A_1.north east) to (P_2);

\draw[-latex,thick,black] (t0) to (t2);

\end{tikzpicture}
\caption{Dynamic decision network encoding the asset-twin coupled dynamical system. Circle nodes denote random variables, square nodes denote actions, and diamond nodes denote the objective function. Bold outlines denote observed quantities, while thin outlines denote estimated quantities. Directed solid edges represent the variables' dependencies encoded via conditional probability distributions, while directed dashed edges represent the variables' dependencies encoded via deterministic functions.
\label{fig:graph}}
\end{figure}
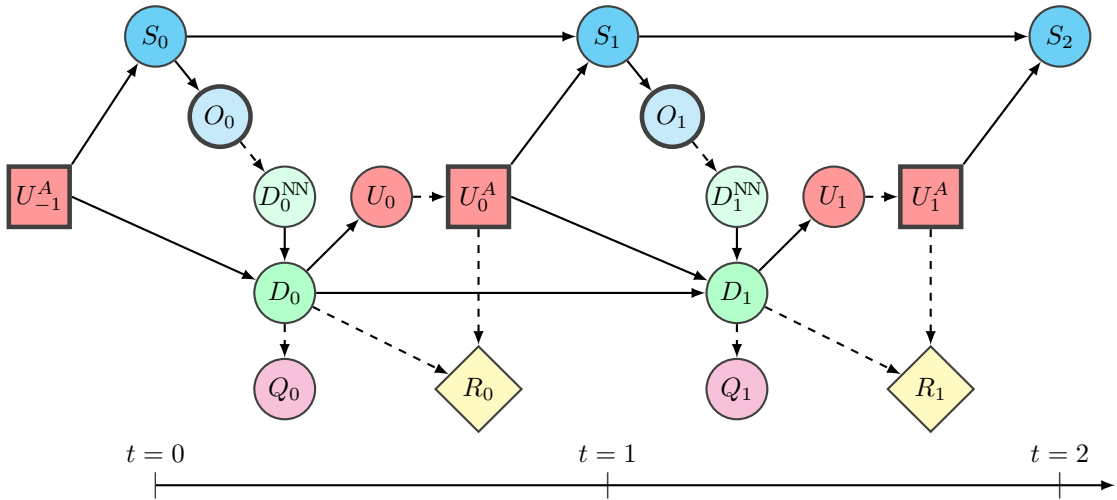

The physical state $S_t\sim p(s_t)$, with $s_t$ denoting the realization of the random variable $S_t$ at time $t$, encapsulates the variability in the health state of the asset, which is usually only partially observable. The probability distribution encoding the relative likelihood that $S_t=s_t$, for any possible $s_t$, is denoted with $p(s_t)$. The digital state $D_t\sim p(d_t)$ is characterized by those parameters employed to capture the variability of the physical asset by means of the computational models comprising the DT. In our framework, the digital state is given as a vector of length two, describing the presence/location and magnitude of damage in the asset. The \textit{physical-to-digital} information flow is governed by the observed data $O_t=o_t$, which are assimilated by the DT to update the digital state. The assimilation is carried out using the DL models described in \sez\ref{sec:DL_models}, providing a first estimate of the digital state $D^\text{N\hspace{-1px}N}_t\sim p(d^\text{N\hspace{-1px}N}_t)$. This estimation is then used in a Bayesian inference formulation, together with the prior belief $D_{t-1}$ from the previous time step, to estimate an updated digital state $D_{t}$ according to a control-dependent transition dynamics model describing how the digital state is expected to evolve. The updated digital state can thus be exploited to compute quantities of interest $Q_t\sim p(q_t)$, such as modal quantities or other response features, through the computational models comprising the DT. For instance, quantities of interest can be useful to perform posterior predictive checks on the tracking capabilities of the DT to assess how it matches the reality, by comparing sensor measurements with the corresponding posterior estimates. However, we point out that this capability is not exploited in the present work, and that the $Q_t$ node is kept in the graph in agreement with the foundational model proposed in~\cite{art:kapteyn2021probabilistic}. Nevertheless, the updated digital state $D_t$ is eventually exploited to inform the \textit{digital-to-physical} information flow, in the form of control inputs; in \fig\ref{fig:graph}, $U_t\sim p(u_t)$ and $U^A_t=u^A_t$ denote the belief about what action to take and the control input effectively enacted on the asset, respectively. At each time step, $U_t$ is estimated according to a health-dependent control policy, that maps the belief over the digital state onto the control actions feeding back to the physical asset. Finally, the reward $R_t\sim p(r_t)$ quantifies the performance of the asset for the time step and can be equivalently perceived as a negative cost to be maximized.

The key assumptions behind our PGM are that the physical state is only observable indirectly via the sensed structural response, and the physical and digital states evolve according to a Markovian process. This implies that the conditional probabilities associated with the random variables at one time step depend only on the random variables at the previous time step, and are independent of all past states. The resulting graph topology encodes a conditional independence structure that allows us to conveniently cast the asset tracking within a sequential Bayesian inference framework. Indeed, by exploiting conditional independence and Bayes rule, the joint distributions over variables can be factorized up to the current time step $t_c$, as follows:
\begin{align}
&p(D^{\text{N\hspace{-1px}N}}_0,\dots,D^{\text{N\hspace{-1px}N}}_{t_c},D_0,\dots,D_{t_c},Q_0,\dots,Q_{t_c},R_0,\dots,R_{t_c},U_{0},\dots,U_{t_c}|o_0,\dots,o_{t_c},u^A_0,\dots,u^A_{t_c}) \nonumber \\  &\propto\prod_{t=0}^{t_c}\Bigl[\phi_t^\text{data}\phi_t^\text{history}\phi_t^\text{N\hspace{-1px}N}\phi_t^\text{QoI}\phi_t^\text{control}\phi_t^\text{reward} \Bigr],
\label{eq:factorization}
\end{align}
with factors:
\begin{align}
&\phi_t^\text{data} = p(O_{t} = o_{t} | D^\text{N\hspace{-1px}N}_t),\\
&\phi_t^\text{history} = p(D_t | D_{t-1}, U^A_{t-1} = u^A_{t-1}),\\
&\phi_t^\text{N\hspace{-1px}N} = p(D_t|D^{\text{N\hspace{-1px}N}}_{t}),\\
&\phi_t^\text{QoI} = p(Q_t|D_t),\\
&\phi_t^\text{control}=p(U_t|D_t),\\
&\phi_t^\text{reward} = p(R_t | D_t,U^A_t =u^A_t).
\end{align}
Herein, $\phi_t^\text{data}$ encodes the assimilation of observed data through the DL models underlying the identification of the structural health. $\phi_t^\text{history}$ and $\phi_t^\text{N\hspace{-1px}N}$ factorize the belief about the digital state $D_t$, conditioned on the digital state at the previous time step $D_{t-1}$, the last enacted control input $U^A_{t-1} =u^A_{t-1}$, and the data assimilation outcome $D^{\text{N\hspace{-1px}N}}_{t}$. In our PGM, the spaces of the digital states and control inputs are discrete, thus the relevant causal relationships are modeled by means of conditional probability tables (CPTs). In particular, $\phi_t^\text{history}$ plays the role of a predictor forward in time, parametrized by means of a control-dependent CPT describing how the digital state is expected to evolve. Such a CPT should embody any a priori knowledge that the DT designer has with respect to the asset and the relevant operational conditions. $\phi_t^\text{history}$ can be estimated offline from historical data, see e.g.,~\cite{art:state_trans_1,art:state_trans_2}, or learned online from experience. On the other hand, $\phi_t^\text{N\hspace{-1px}N}$ updates the digital state estimate to account for data assimilation. This is encoded by means of a CPT mapping the estimate $D^{\text{N\hspace{-1px}N}}_{t}$ provided by the DL models, onto a belief about $D_t$. Such a CPT is a confusion matrix measuring the offline (expected) performance of the DL models in correctly identifying the digital state among all the possible outcomes of $D_t$. $\phi_t^\text{QoI}$ and $\phi_t^\text{reward}$ respectively encapsulate the evaluation of the computational models comprising the DT to estimate quantities of interest, and the computation of the reward function quantifying the performance of the asset. Finally, the control factor $\phi_t^\text{control}$ is encoded by means of a health-dependent control policy $\pi(D_t)$ computed as described in the following. Since the spaces of the unobserved variables are discrete, we can propagate and update the relative belief exactly with a single pass of the sum-product message-passing algorithm \cite{book:2009probabilistic}.

The control policy $\pi(D_t)$ is computed offline under the simplifying assumption of sufficient sensing capability to provide an accurate estimate of the structural health, allowing us to decouple the sensing and control problems. This involves solving the planning problem induced by the expected evolution of the structural health, maximizing the expected reward over the planning horizon. Considering an infinite planning horizon, this can be stated as the optimization problem:
\begin{equation}
\pi(D_t) = \argmax_\pi{\sum_{t=0}^{+\infty}\gamma^t\mathbb{E}[R_t]},
\label{eq:optimization}
\end{equation}
where $\gamma\in[0,1]$ is the discount factor. Here, this is solved using the dynamic-programming value iteration algorithm~\cite{book:Sutton}. The reward function to be optimized is chosen as:
\begin{equation}
R_t(U_t,D_t) = R_t^\text{control}(U_t) + \alpha R_t^\text{health}(D_t).
\label{eq:reward}
\end{equation}
Herein, $R_t^\text{control}(U_t)$ and $R_t^\text{health}(D_t)$ quantify the rewards relative to control inputs and health state, respectively, and $\alpha\in\mathbb{R}$ is a weighting factor, useful to tune the trade-off between risk-averse and risk-seeking behavior. After learning $\pi(D_t)$, $U^A_t$ is selected as the best point estimate of $U_t$.

\begin{figure}[h!]
\center
\begin{tikzpicture}[scale=.85, every node/.style={scale=.9}]

\node [D_node] (D_0) at (0,0) {};
\node [] at (0,0) {$D_{0}$};
\node [U_node_prob] (U_0) at (1.5,1.5) {};
\node [] at (1.5,1.5) {$U_{0}$};
\node [Q_node] (Q_0) at (0,-1.5) {};
\node [] at (0,-1.5) {$Q_{0}$};
\node [R_node] (R_0) at (1.5,-1.5) {};
\node [] at (1.5,-1.5) {$R_{0}$};

\node [D_node] (D_1) at (3,0) {};
\node [] at (3,0) {$D_{1}$};
\node [U_node_prob] (U_1) at (4.5,1.5) {};
\node [] at (4.5,1.5) {$U_{1}$};
\node [Q_node] (Q_1) at (3,-1.5) {};
\node [] at (3,-1.5) {$Q_{1}$};
\node [R_node] (R_1) at (4.5,-1.5) {};
\node [] at (4.5,-1.5) {$R_{1}$};

\node [D_node] (D_2) at (6,0) {};
\node [] at (6,0) {$D_{2}$};
\node [U_node_prob] (U_2) at (7.5,1.5) {};
\node [] at (7.5,1.5) {$U_{2}$};
\node [Q_node] (Q_2) at (6,-1.5) {};
\node [] at (6,-1.5) {$Q_{2}$};
\node [R_node] (R_2) at (7.5,-1.5) {};
\node [] at (7.5,-1.5) {$R_{2}$};

\node [D_node] (D_3) at (9,0) {};
\node [] at (9,0) {$D_{3}$};
\node [U_node_prob] (U_3) at (10.5,1.5) {};
\node [] at (10.5,1.5) {$U_{3}$};
\node [Q_node] (Q_3) at (9,-1.5) {};
\node [] at (9,-1.5) {$Q_{3}$};
\node [R_node] (R_3) at (10.5,-1.5) {};
\node [] at (10.5,-1.5) {$R_{3}$};

\node [] (t0) at (-.14,-3) {};
\node [] () at (0,-3) {$|$};
\node [] () at (0,-2.5) {$t=t_c$};
\node [] () at (3,-3) {$|$};
\node [] () at (3,-2.5) {$t=t_c+1$};
\node [] () at (6,-3) {$|$};
\node [] () at (6,-2.5) {$t=t_c+2$};
\node [] () at (9,-3) {$|$};
\node [] () at (9,-2.5) {$t=t_c+3$};
\node [] (t2) at (13,-3) {};

\draw[-latex,thick,black] (D_0) to (U_0);
\draw[-latex,thick,black] (U_0) to (D_1);
\draw[-latex,thick,black] (D_0) to (D_1);
\draw[-latex,dashed,thick,black] (D_0) to (Q_0);
\draw[-latex,dashed,thick,black] (D_0) to (R_0);
\draw[-latex,dashed,thick,black] (U_0) to (R_0);

\draw[-latex,thick,black] (D_1) to (U_1);
\draw[-latex,thick,black] (U_1) to (D_2);
\draw[-latex,thick,black] (D_1) to (D_2);
\draw[-latex,dashed,thick,black] (D_1) to (Q_1);
\draw[-latex,dashed,thick,black] (D_1) to (R_1);
\draw[-latex,dashed,thick,black] (U_1) to (R_1);

\draw[-latex,thick,black] (D_2) to (U_2);
\draw[-latex,thick,black] (U_2) to (D_3);
\draw[-latex,thick,black] (D_2) to (D_3);
\draw[-latex,dashed,thick,black] (D_2) to (Q_2);
\draw[-latex,dashed,thick,black] (D_2) to (R_2);
\draw[-latex,dashed,thick,black] (U_2) to (R_2);

\draw[-latex,thick,black] (D_3) to (U_3);
\draw[-latex,dashed,thick,black] (D_3) to (Q_3);
\draw[-latex,dashed,thick,black] (D_3) to (R_3);
\draw[-latex,dashed,thick,black] (U_3) to (R_3);

\draw[-latex,thick,black] (t0) to (t2);

\end{tikzpicture}
\caption{Dynamic decision network employed to predict the future evolution of the digital state and the associated uncertainty. Circle nodes denote random variables, and diamond nodes denote the objective function. Directed solid edges represent the variables' dependencies encoded via conditional probability distributions, while directed dashed edges represent the dependencies encoded via deterministic functions.\label{fig:graph_prediction}}
\end{figure}
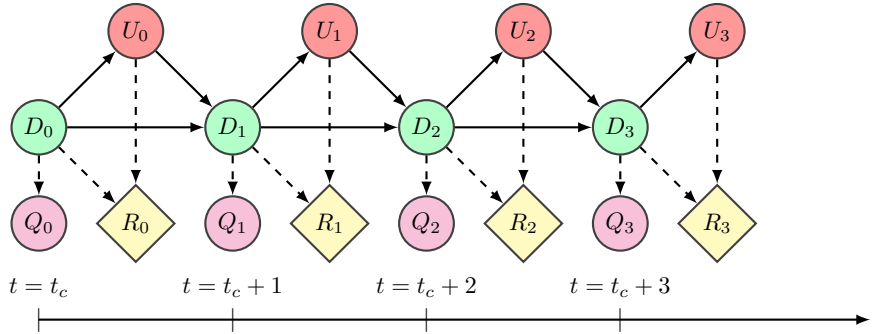

Starting from the updated digital state $D_{t_c}$ at the current time step $t_c$, future prediction is achieved by unrolling until a prediction time $t_p$ the portion of PGM relative to $D_t,Q_t,R_t$, and $U_t$ (see \fig\ref{fig:graph_prediction}). All other nodes are removed from the prediction graph, as neither data assimilation nor actions are performed on the asset while forecasting its evolution. The factorization in \eq\eqref{eq:factorization} can be extended over the prediction horizon as:
\begin{align}
&p(D^{\text{N\hspace{-1px}N}}_0,\ldots,D^{\text{N\hspace{-1px}N}}_{t_c},D_0,\ldots,D_{t_p},Q_0,\ldots,Q_{t_p},R_0,\ldots,R_{t_p},U_{0},\ldots,U_{t_p}|o_0,\ldots,o_{t_c},u^A_0,\ldots,u^A_{t_c}) \nonumber \\
&\propto\prod_{t=0}^{t_p}\Bigl[\phi_t^\text{history}\phi_t^\text{QoI}\phi_t^\text{control}\phi_t^\text{reward} \Bigr]\prod_{t=0}^{t_c}\Bigl[\phi_t^\text{data}\phi_t^\text{N\hspace{-1px}N}\Bigr].
\label{eq:factorization_pred}
\end{align}

The algorithmic description of the online phase of the proposed digital twinning framework is reported in \alg\ref{al:online}. The operations repeat each time new observational data are provided. Note that the considered PGM digital twinning framework is general, and can easily be adapted to deal with physical assets other than civil engineering structures by reorganizing the topology of the graph, if necessary.

\begin{algorithm}[ht]
\hspace*{\algorithmicindent} \textbf{Input}: observational data $O_t=o_t$
\begin{algorithmic}[1]
\State assimilate $o_t$ with the DL models to provide $D^{\text{N\hspace{-1px}N}}_{t}=d^{\text{N\hspace{-1px}N}}_{t}$.\Comment{$(O_t)\rightarrow(D^{\text{N\hspace{-1px}N}}_{t})$}
\State infer $D_{t}$ and $U_{t}$ by updating $d_{t-1}$, given $u^A_{t-1}$, $d^{\text{N\hspace{-1px}N}}_{t}$, and the CPTs encoding $\phi_t^\text{history}$, $\phi_t^\text{N\hspace{-1px}N}$ and $\phi_t^\text{control}$.\Comment{$(D_{t-1},D^{\text{N\hspace{-1px}N}}_{t},U^A_{t-1},)\rightarrow(D_{t},U_{t})$}
\State infer the future evolution of $D_{t}$ and $U_{t}$, given the updated $d_{t}$, and the CPTs encoding $\phi_t^\text{history}$ and $\phi_t^\text{control}$.\Comment{$(D_{t_c} )\rightarrow(D_{t_p},U_{t_p})$}
\State select $U^A_t=u^A_t$ as the best point estimate of $U_{t}=u_{t}$.\Comment{$(U_{t})\rightarrow(U^A_{t})$}
\State \Return control input to be enacted $u^A_t$, expected evolution of $D_{t}$ and $U_{t}$.
\end{algorithmic}
\caption{Online phase -- algorithmic description}
\label{al:online}
\end{algorithm}

\subsection{Numerical models for simulation-based damage identification}
\label{sec:numerical_models}

As anticipated in the previous section, the assimilation of structural response data to identify the structural state is carried out through DL models. A simulation-based strategy is exploited to train DL models on the basis of vibration responses. The training data are numerically generated by simulating physics-based models so that the effect of damage on the structural response can be systematically reproduced~\cite{art:Taddei_Patera}. In particular, the structure to be monitored is modeled as a linear-elastic continuum, discretized in space through finite elements. Its dynamic response to the applied loadings, under the assumption of linearized kinematics, is described by the following semi-discretized form of the elasto-dynamic problem:
\begin{equation}
\left\{
\begin{array}{ll}
\mathbf{M}\ddot{\mathbf{x}}(t) + \mathbf{C}(\boldsymbol{\mu})\dot{\mathbf{x}}(t) + \mathbf{K}(\boldsymbol{\mu})\mathbf{x}(t)=\mathbf{f}(t,\boldsymbol{\mu}), & t\in(0,T)\\ \mathbf{x}(0)=\mathbf{x}_{0}, & \\ \dot{\mathbf{x}}(0)=\dot{\mathbf{x}}_{0}, & 
\end{array}
\right.
\label{eq:HF_model}
\end{equation}
which is referred to as the full-order model (FOM). Here $t\in(0,T)$ denotes time; $\mathbf{x}(t),\dot{\mathbf{x}}(t),\ddot{\mathbf{x}}(t)\in\mathbb{R}^{N_\text{FE}}$ are the vectors of nodal displacements, velocities and accelerations, respectively; $N_\text{FE}$ is the number of degrees of freedom (dofs); $\mathbf{M}\in\mathbb{R}^{N_\text{FE}\times N_\text{FE}}$ is the mass matrix; $\mathbf{C}(\boldsymbol{\mu})\in\mathbb{R}^{N_\text{FE}\times N_\text{FE}}$ is the damping matrix, assembled according to the Rayleigh's model; $\mathbf{K}(\boldsymbol{\mu})\in\mathbb{R}^{N_\text{FE}\times N_\text{FE}}$ is the stiffness matrix; $\mathbf{f}(t,\boldsymbol{\mu})\in\mathbb{R}^{N_\text{FE}}$ is the vector of nodal forces induced by the external loadings; and $\mathbf{x}_{0}$ and $\dot{\mathbf{x}}_{0}$ are the initial conditions (at $t=0$), in terms of nodal displacements and velocities, respectively. The mass matrix $\mathbf{M}$ is not a function of $\boldsymbol{\mu}$ because the mass properties of the structure are unaffected by the employed damage description or by the operational conditions. The solution of \pb\eqref{eq:HF_model} is advanced in time using the Newmark integration scheme (constant average acceleration method)~\cite{art:Newmark}, to provide $\mathbf{x}_l$, $\dot{\mathbf{x}}_l$ and $\ddot{\mathbf{x}}_l$, for $l=1,\ldots,L$, with $\mathbf{x}_l$ being the vector of nodal displacements at time $l$.

With reference to civil structures, we focus on the early detection of damage patterns characterized by a small evolution rate, whose prompt identification can reduce lifecycle costs and increase the safety and availability of the structure. In this context, a localized reduction of the material stiffness stands as the simplest damage mechanism resulting from a time scale separation between damage growth and damage assessment, see e.g.,~\cite{art:TEUGHELS2002,book:ML_perspective,art:PANDEY19943}. Here, local stiffness reduction is obtained by parametrizing the stiffness matrix via two variables $y\in\mathbb{N}$ and $\delta\in\mathbb{R}$, included in the parameter vector $\boldsymbol{\mu}$, respectively describing the location and magnitude of the applied stiffness reduction, similarly to~\cite{art:AMSES20,art:Metodologico,art:Torzoni_temperature}. In particular, $y\in\lbrace0,\ldots,N_y\rbrace$ labels the specific damage region, among a set of predefined $N_y$ damage locations, where $y=0$ identifies the damage-free baseline. The parameter $\delta\in\mathbb{R}$ describes the magnitude of the stiffness reduction taking place within the predesignated region associated with $y$.

As $N_\text{FE}$ increases, the computational cost associated with the solution of the FOM for any sampled $\boldsymbol{\mu}$ also grows, and the generation of synthetic datasets becomes prohibitive. To address this challenge, a projection-based reduced-order model (ROM) is exploited in place of the FOM to speed up the offline dataset population phase, similarly to~\cite{art:Metodologico,art:Torzoni_temperature}. The ROM is obtained by a proper orthogonal decomposition (POD)-Galerkin reduced basis method~\cite{book:RB, chinestaenc2017, book:degruyter1, book:RozzaStabileBallarin2022}. This reduced-order modeling strategy is chosen because POD has been investigated and validated in the context of structural dynamics~\cite{art:Kerschen_1,art:Kerschen_2} and structural analysis~\cite{art:guo2018reduced, art:tezzele2021multi}, its appealing offline-online decoupling, and the availability of efficient criteria for the selection of POD basis functions. It is worth noting that alternative reduced-order modeling approaches can also be employed to alleviate the computational burden during the offline dataset generation. For instance, one could use spectral POD~\cite{sieber2016spectral, towne2018spectral, lario2022neural}, or Grassmannian diffusion maps~\cite{dos2022grassmannian}, as viable alternatives to the reduced basis method.

The ROM approximation to the solution of \pb\eqref{eq:HF_model} is obtained by linearly combining $N_\text{RB}\ll N_\text{FE}$ POD basis functions $\mathbf{w}_{k}\in\mathbb{R}^{N_\text{FE}}$, \mbox{$k=1,\ldots,N_\text{RB}$}, as $\mathbf{x}(t,\boldsymbol{\mu})\approx\mathbf{W}\widehat{\mathbf{x}}(t,\boldsymbol{\mu})$, where $\mathbf{W}=[\mathbf{w}_1,\ldots,\mathbf{w}_{N_\text{RB}}] \in \mathbb{R}^{N_\text{FE} \times N_\text{RB}}$ is the basis matrix collecting the POD basis functions and $\widehat{\mathbf{x}}(t,\boldsymbol{\mu})\in \mathbb{R}^{N_\text{RB}}$ is the vector of unknown POD coefficients. By enforcing the orthogonality between the residual and the subspace spanned by the first $N_\text{RB}$ POD modes through a Galerkin projection, the following $N_\text{RB}$-dimensional semi-discretized form is obtained:
\begin{equation}
\left\{
\begin{array}{ll}
\mathbf{M}_{r}\ddot{\widehat{\mathbf{x}}}(t) + \mathbf{C}_r(\boldsymbol{\mu})\dot{\widehat{\mathbf{x}}}(t) + \mathbf{K}_{r}(\boldsymbol{\mu})\widehat{\mathbf{x}}(t)=\mathbf{f}_{r}(t,\boldsymbol{\mu}), &  t\in(0,T)\\ \widehat{\mathbf{x}}(0)=\mathbf{W}^\top\mathbf{x}_{0}, & \\ \dot{\widehat{\mathbf{x}}}(0)=\mathbf{W}^\top\dot{\mathbf{x}}_{0}. &
\end{array}
\right.
\label{eq:LF_model}
\end{equation}
The solution of this reduced-order system is advanced in time using the same strategy employed for the FOM model, and then projected onto the original FOM space as $\mathbf{x}(t,\boldsymbol{\mu})\approx\mathbf{W}\widehat{\mathbf{x}}(t,\boldsymbol{\mu})$. Here, reduced matrices $\mathbf{M}_{r}$, $\mathbf{C}_r$, and $\mathbf{K}_{r}$, and the reduced vector $\mathbf{f}_{r}$ play the same role as their high-fidelity counterparts, yet with dimension $N_\text{RB} \times N_\text{RB}$ instead of $N_\text{FE} \times N_\text{FE}$, according to the following relationships:
\begin{equation}
\begin{array}{lll}
\mathbf{M}_{r} \equiv \mathbf{W}^\top \mathbf{M} \mathbf{W}, &\quad &\mathbf{C}_r(\boldsymbol{\mu}) \equiv \mathbf{W}^\top \mathbf{C}(\boldsymbol{\mu})\mathbf{W},\\ 
\mathbf{K}_{r}(\boldsymbol{\mu}) \equiv \mathbf{W}^\top \mathbf{K}(\boldsymbol{\mu})\mathbf{W}, &\quad&
\mathbf{f}_{r}(t,\boldsymbol{\mu}) \equiv \mathbf{W}^\top \mathbf{f}(t,\boldsymbol{\mu}).
\end{array}
\end{equation}

The basis matrix $\mathbf{W}$ is obtained by POD, exploiting the so-called method of snapshots as follows. First, a snapshot matrix $\mathbf{S}=[\mathbf{x}_1,\ldots,\mathbf{x}_{N_\text{S}}]\in\mathbb{R}^{N_\text{FE}\times N_\text{S}}$ is assembled from $N_\text{S}$ solution snapshots, computed by integrating in time the FOM solution for different values of parameters $\boldsymbol{\mu}$. The computation of an optimal reduced basis is then carried out by factorizing $\mathbf{S}$ through a singular value decomposition. We use a standard energy-based criterion to set the order $N_\text{RB}$ of the approximation. For further details see, e.g.,~\cite{book:RB,book:RozzaStabileBallarin2022,art:Torzoni_MF,art:Torzoni_DML}.

To populate the training dataset $\mathcal{D}$, the parametric space of vector $\boldsymbol{\mu}$ is taken as uniformly distributed, and then sampled via the Latin hypercube rule. The number of samples is equal to the number $I$ of instances collected in $\mathcal{D}$ as: 
\begin{equation}
\mathcal{D}=\lbrace(\mathbf{U}_i,y_i,\delta_i)\rbrace_{i=1}^{I},
\label{eq:Dataset}
\end{equation}
where the vibration recordings $\mathbf{U}_i$ associated with the $i$-th sampling of $\boldsymbol{\mu}$, with $i=1,\ldots,I$, are labeled by the corresponding values of $y_i$ and $\delta_i$, and are obtained as follows. With reference to displacement recordings, nodal values in $(0,T)$ are first collected as $\mathbf{V}_i=[\mathbf{W}\widehat{\mathbf{x}}_1,\ldots,\mathbf{W}\widehat{\mathbf{x}}_{L}]_i \in \mathbb{R}^{N_\text{FE} \times L}$ by solving \pb\eqref{eq:LF_model}. The relevant vibration recordings $\mathbf{U}_i$ are then obtained as:
\begin{equation}
\qquad \mathbf{U}_i=(\mathbf{T}\mathbf{V}_i)^\top,
\label{eq:bool_sensors}
\end{equation}
where $\mathbf{T}\in\mathbb{R}^{N_u \times N_\text{FE}}$ is a Boolean matrix whose $(n,m)$--th entry is equal to $1$ only if the $n$--th sensor output coincides with the $m$--th dof. In order to mimic the measurement noise, each vibration recording in $\mathcal{D}$ is corrupted by adding an independent, identically distributed Gaussian noise, whose statistical properties depend on the target accuracy of the sensors. In the following, the index $i$ will be dropped for ease of notation, unless necessary.

\subsection{Data assimilation via artificial neural networks}
\label{sec:DL_models}

The $\phi_t^\text{data}$ factor in our PGM encodes the assimilation of observed data through the DL models underlying the identification of the structural health. In this section, we describe the adopted DL models, the aspects related to their training, and how they are used to assimilate observational data to detect, locate, and quantify the presence of structural damage.

Every time new observational data $\mathbf{U}$ are acquired, they are first processed with a classification model $\text{N\hspace{-1px}N}_\text{CL}$ to address damage detection/localization. Classification involves the prediction of an output class to categorize a given input. Here, the classes are those described through the $y$ parameter. Whenever a damage is identified in the $j$-th region, $j=1,\ldots,N_y$, the observational data $\mathbf{U}$ are further processed with regression models $\text{N\hspace{-1px}N}^j_\text{RG}$, one for each damageable region, to quantify the associated amount of damage $\delta$.

The aforementioned classification and regression tasks are addressed by means of DL models. The use of DL models for SHM purposes has the advantage of automating the feature engineering stage characterizing the pattern recognition paradigm for SHM~\cite{book:ML_perspective,book:Bishop}. Indeed, a DL model is trained to select and extract optimized damage-sensitive features from raw sensor recordings through an end-to-end learning process. Moreover, since the DL model is learned offline, the structural state can be next assessed in real-time regardless of considering continuous or discrete variables, which would be difficult to achieve with other optimization techniques, such as nonlinear programming, stochastic optimization, and metaheuristic methods.

The model $\text{N\hspace{-1px}N}_\text{CL}$ addresses the multi-class classification task underlying the damage detection/localization problem, namely $\text{N\hspace{-1px}N}_\text{CL}:\mathbf{U}\rightarrow\boldsymbol{b}\in\mathbb{R}^{N_y+1}$. The target label $\boldsymbol{b}$ categorizes one of the $N_y+1$ predefined damage scenarios described through parameter $y$. In particular, $\boldsymbol{b}$ is a one-hot encoding $\boldsymbol{b} = [ b^0,\ldots,b^{N_y}]^\top$, with entries $b^m$ equal to $1$ if the target class $y$ is $m$ and $0$ otherwise, with $m=0,\ldots,N_y$. This is needed because DL models cannot operate on nominal data directly. They require all input variables and output variables to be numeric. The one-hot encoding converts the nominal feature described by the $y$ parameter into a multidimensional binary vector. The number of dimensions corresponds to the number of categories, and each category gets its dimension. Each category is encoded by mapping it to a vector in which the entry corresponding to the category’s dimension is $1$, and the rest are $0$. 

The estimated counterpart of $\boldsymbol{b}$ is obtained as $\widehat{\boldsymbol{b}}=\text{N\hspace{-1px}N}_\text{CL}(\mathbf{U})$. By employing a Softmax activation function for the output layer of $\text{N\hspace{-1px}N}_\text{CL}$, the entries of $\widehat{\boldsymbol{b}} = ( \widehat{b}^0,\ldots,\widehat{b}^{N_y})^\top\in\mathbb{R}^{N_y+1}$ are interpreted as the confidence levels $\widehat{b}^m$ by which $\mathbf{U}$ is assigned to the $m$-th damage class, with $m=0,\ldots,N_y$. In particular, the Softmax activation function converts the real-valued vector $\boldsymbol{a}=(a^0,\ldots,a^{N_y})^\top\in\mathbb{R}^{N_y+1}$, provided by the output layer of $\text{N\hspace{-1px}N}_\text{CL}$, into a discrete probability distribution as:
\begin{equation}
\widehat{\boldsymbol{b}}=\text{Softmax}(\boldsymbol{a}), \qquad \text{with}\quad\widehat{b}^m(\boldsymbol{a}) = \frac{\text{exp}(a^m)}{\sum_{k=0}^{N_y} \text{exp}(a^{k})}, \quad m=0,\ldots,N_y.
\label{eq:softmax}
\end{equation}
 When $\text{N\hspace{-1px}N}_\text{CL}$ is exploited for prediction, the most likely class is selected as the one that best categorizes the processed measurements $\mathbf{U}$.

The model $\text{N\hspace{-1px}N}^j_\text{RG}$ addresses the regression task underlying the damage quantification problem, namely $\text{N\hspace{-1px}N}^j_\text{RG}:\mathbf{U}\rightarrow\delta\in\mathbb{R}$, with $j=1,\ldots,N_y$. The estimated counterpart of $\delta$ is obtained as $\widehat{\delta}=\text{N\hspace{-1px}N}^j_\text{RG}(\mathbf{U})$. Hence, the regression models, one for each damageable region, map the vibration recordings $\mathbf{U}$ associated with the $j$-th damage region, onto the estimated magnitude of the stiffness reduction taking place within the relative damage region. Since all $\text{N\hspace{-1px}N}^j_\text{RG}$ models are learned following the same procedure, the index $j$ will be dropped in the following for ease of notation. 

Since the space of digital states in the PGM is discrete, the outcomes of $\text{N\hspace{-1px}N}_\text{CL}$ and $\text{N\hspace{-1px}N}_\text{RG}$ are accommodated within the PGM by discretizing the range in which the damage level $\delta$ can take values in $N_\delta$ uniform intervals, thus resulting in $N_d=1+N_\delta N_y$ possible states. The same reasoning is followed to compute the confusion matrix encoding the $\phi_t^\text{N\hspace{-1px}N}$ factor. In particular, $\phi_t^\text{N\hspace{-1px}N}$ measures the offline performance of $\text{N\hspace{-1px}N}_\text{CL}$ and $\text{N\hspace{-1px}N}_\text{RG}$ in assimilating noisy FOM data to classify the digital state, among the $N_d$ possible outcomes of $D_t$.

The models $\text{N\hspace{-1px}N}_\text{CL}$ and $\text{N\hspace{-1px}N}_\text{RG}$ are trained separately. The datasets dedicated to the training of $\text{N\hspace{-1px}N}_\text{CL}$ and $\text{N\hspace{-1px}N}_\text{RG}$ are derived from dataset $\mathcal{D}$ in \eq\eqref{eq:Dataset} as follows. The dataset used to learn $\text{N\hspace{-1px}N}_\text{CL}$ is obtained from \eq\eqref{eq:Dataset}, as
\begin{equation}
    \mathcal{D}_{\text{CL}}=\lbrace (\mathbf{U}_i,\boldsymbol{b}_i)\rbrace^I_{i=1}.
    \label{eq:Dataset_cl}
\end{equation}
The dataset used to learn $\text{N\hspace{-1px}N}_\text{RG}$ is derived from \eq\eqref{eq:Dataset}, as
\begin{equation}
\mathcal{D}_{\text{RG}}=\lbrace (\mathbf{U}_{i_{\text{RG}}},\delta_{i_{\text{RG}}})\rbrace^{I_{\text{RG}}}_{i_{\text{RG}}=1},
    \label{eq:Dataset_rg}
\end{equation}
where $I_{\text{RG}}$ is the number of training instances in $\mathcal{D}_{\text{RG}}$, all characterized by a structural damage within the same predefined region.

The set of weights and biases parametrizing $\text{N\hspace{-1px}N}_\text{CL}$ is denoted as $\boldsymbol{\Theta}_\text{CL}$. This is optimized minimizing the probabilistic categorical cross-entropy~\cite{book:DL_book,art:AMSES20} $\mathcal{L}_\text{CL}$ between the predicted and target class labels over $\mathcal{D}_\text{CL}$:
\begin{equation}
\mathcal{L}_\text{CL}(\boldsymbol{\Theta}_\text{CL},\mathcal{D}_\text{CL}) = - \frac{1}{I}\sum_{i=1}^{I} \sum_{m=0}^{N_y} b^m_i\ \text{log} (\widehat{b}^m_i),
\label{eq:cross_entropy}
\end{equation}
which provides a measure of the distance between the discrete probability distribution describing $\boldsymbol{b}$, and its estimated counterpart $\widehat{\boldsymbol{b}}=\text{N\hspace{-1px}N}_\text{CL}(\mathbf{U})$.

The set of weights and biases $\boldsymbol{\Theta}_\text{RG}$ parametrizing $\text{N\hspace{-1px}N}_\text{RG}$ is learned through the minimization of the following mean squared error loss function:
\begin{equation}
\mathcal{L}_\text{RG}(\boldsymbol{\Theta}_\text{RG},\mathcal{D}_\text{RG}) =\frac{1}{I_{\text{RG}}}\sum^{I_{\text{RG}}}_{i_{\text{RG}}=1}(\delta_{i_{\text{RG}}} - \widehat{\delta}_{i_{\text{RG}}})^2,
\label{eq:mse_loss}
\end{equation}
which provides a measure of the distance between the target magnitude of the stiffness reduction $\delta$, and its approximated counterpart $\widehat{\delta}=\text{N\hspace{-1px}N}_\text{RG}(\mathbf{U})$. 

The algorithmic description of the procedures and computations characterizing the preliminary offline phase of the proposed digital twinning framework is reported in \alg\ref{al:offline}. The implementation details of the deep learning models are reported in Appendix~\ref{sec:implementation}.

\begin{algorithm}[ht]
\hspace*{\algorithmicindent} \textbf{Input}: parametrization of the operational and damage conditions\\
\hspace*{50pt} PGM implementing the prediction graph
\begin{algorithmic}[1]
\State set up the physics-based numerical model of the structure to be monitored.
\State assemble the snapshot matrix of the structural response via FOM analyses.
\State compute the POD basis functions via singular value decomposition of the snapshots matrix.
\State use the ROM to populate the training dataset $\mathcal{D}$ with vibration recordings at sensor location.
\State use the recordings and labels in $\mathcal{D}$ to derive $\mathcal{D}_\text{CL}$ and $\mathcal{D}_\text{RG}$.
\State train the classification model $\text{N\hspace{-1px}N}_\text{CL}$ on $\mathcal{D}_\text{CL}$ and the regression models $\text{N\hspace{-1px}N}_\text{RG}$ on $\mathcal{D}_\text{RG}$.
\State test the generalization capabilities of $\text{N\hspace{-1px}N}_\text{CL}$ and $\text{N\hspace{-1px}N}_\text{RG}$ on noisy FOM data.
\State compute the confusion matrix encoding the $\phi_t^\text{N\hspace{-1px}N}$ factor.
\State compute the control policy $\pi(D_t)$ by solving the planning problem induced by the PGM.
\State \Return trained DL models, $\phi_t^\text{N\hspace{-1px}N}$ factor, control policy $\pi(D_t)$.
\end{algorithmic}
\caption{Preliminary offline phase -- algorithmic description}
\label{al:offline}
\end{algorithm}

\section{Numerical experiments}
\label{sec:experiments}

This section demonstrates the proposed methodology for two test cases: an L-shaped cantilever beam and a railway bridge. 

The FOM and ROM in \pb\eqref{eq:HF_model} and  \pb\eqref{eq:LF_model} are implemented in the \texttt{Matlab} environment, using the \texttt{redbKIT} library~\cite{Redbkit}. The PGM framework for predictive digital twins is implemented in \texttt{Python}, using the \texttt{pgmpy} library~\cite{pgmpy}. All computations have been carried out on a PC featuring an \texttt{AMD Ryzen\textsuperscript{TM} 9 5950X} CPU @ 3.4 GHz and 128 GB RAM. The NN architectures are implemented through the \texttt{Tensorflow}-based \texttt{Keras} API~\cite{chollet2015keras}, and trained on a single \texttt{Nvidia GeForce RTX\textsuperscript{TM} 3080} GPU card.

\subsection{L-shaped cantilever beam}
\label{sec:beam}

\begin{figure}[h!]
\begin{centering}
\includegraphics[width=.74\textwidth]{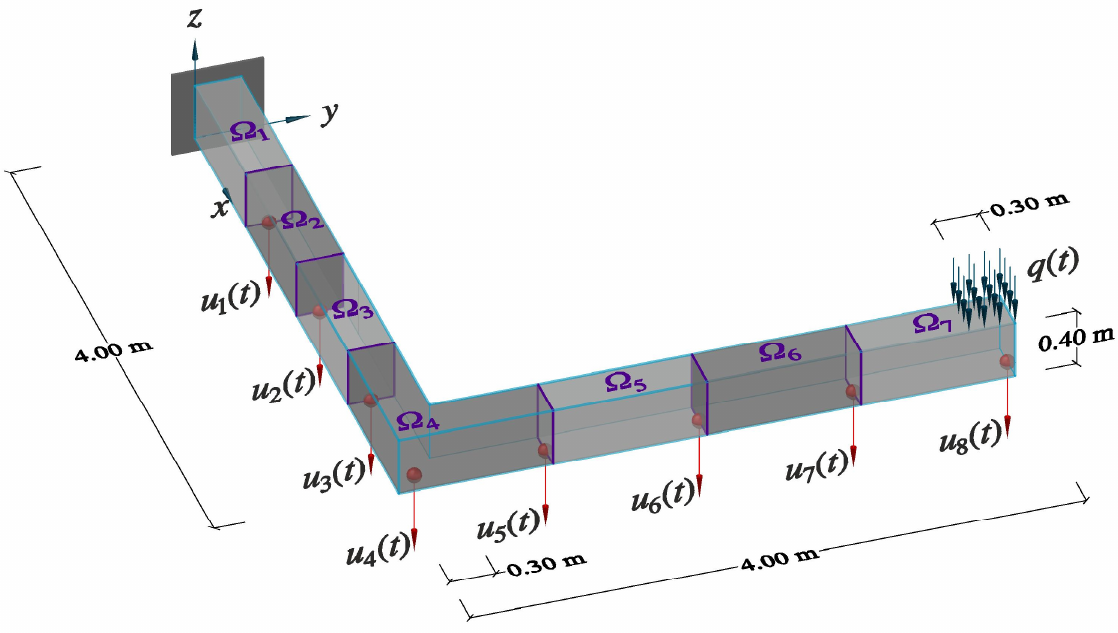}
\caption{L-shaped cantilever beam: details of synthetic recordings related to displacements $u_1(t),\ldots,u_8(t)$, loading condition, and predefined damage regions $\Omega_1,\ldots,\Omega_7$.\label{fig:Maniglia-Model}}
\end{centering}
\end{figure}

The first test case deals with the L-shaped cantilever beam depicted in \fig\ref{fig:Maniglia-Model}. The structure is made of two arms, each one having a length of $4~\textup{m}$, a width of $0.3~\textup{m}$, and a height of $0.4~\textup{m}$. The assumed mechanical properties are those of concrete: Young's modulus $E = 30~\textup{GPa}$, Poisson's ratio $\nu = 0.2$, density $\rho = 2500~\textup{kg/m}^3$. The structure is excited by a distributed vertical load $q(t)$, acting on an area of $(0.3\times0.3)~\textup{m}^2$ close to its tip. The load varies in time according to $q(t) = Q \sin{(2\pi f t)}$, with $Q\in[40, 80]~\textup{kPa}$ and $f\in[10,60]~\textup{Hz}$, respectively being the load amplitude and frequency. Following the setup described in \sez\ref{sec:methodology}, these parameters have a uniform distribution within their respective ranges.

\subsubsection{Dataset assembly}

Synthetic displacement time histories $\mathbf{U}$ are obtained in relation to $N_u=8$ dofs along the bottom surface of the structure, to mimic the monitoring system depicted in \fig\ref{fig:Maniglia-Model}. Each recording is provided for a time interval $(0,T=1~\textup{s})$ with an acquisition frequency $f_\text{s}=200~\textup{Hz}$. Recordings are corrupted with an additive Gaussian noise yielding a signal-to-noise ratio of $100$. 

In addition to the damage-free baseline condition, damage is simulated by considering $N_y=7$ possible damage classes, each referring to a reduction of the material stiffness within a subdomain $\Omega_j$, with $j=1,\ldots,N_y$, as depicted in \fig\ref{fig:Maniglia-Model}. The stiffness reduction can occur with a magnitude $\delta\in[30\%,80\%]$, and is held constant within the considered time interval. 

The FOM is obtained with a finite element discretization using linear tetrahedral elements and resulting in $N_\text{FE}=4659$ dofs. The basis matrix $\mathbf{W}$ is obtained from a snapshot matrix $\mathbf{S}$, assembled through $400$ evaluations of the FOM, at varying values of the input parameters $\boldsymbol{\mu} = (Q,f,y,\delta)^\top$ sampled via Latin hypercube rule. By prescribing a tolerance $\epsilon=10^{-3}$ on the fraction of energy content to be disregarded in the approximation, the order of the ROM approximation turns out to be $N_\text{RB}=56$.

The dataset $\mathcal{D}$ is built with $I=10,000$ instances collected using the ROM. This is then employed to train $\text{N\hspace{-1px}N}_\text{CL}$ and $\text{N\hspace{-1px}N}_\text{RG}$, as described in the previous section. In the absence of experimental data, the testing phase of $\text{N\hspace{-1px}N}_\text{CL}$ and of $\text{N\hspace{-1px}N}_\text{RG}$ is carried out through noise-corrupted FOM solutions. In particular, the asset is monitored by processing batches of $N_\text{obs}=10$ noisy observations, relative to the same damage location $y$ and damage magnitude $\delta$, yet featuring varying operational conditions set by $Q$ and $f$. As the health of the asset evolves over time, the DT assimilates a batch of noisy observations $\lbrace\mathbf{U}_k\rbrace_{k=1}^{N_\text{obs}}$ at each time step, to dynamically estimate the variation in the structural health parameters underlying the digital state.

\subsubsection{Digital twin framework}
The two structural health parameters within the digital state are $\boldsymbol{d}=(y,\delta)^\top$. In order to accommodate the outcome of the DL models within the PGM and to compute the CPT encoding the $\phi_t^\text{N\hspace{-1px}N}$ factor, the range in which the damage level $\delta$ can take values is discretized in $N_\delta=6$ intervals $\lbrace[30\%,35\%],[35\%,45\%],[45\%,55\%],[55\%,65\%],[65\%,75\%],[75\%,80\%]\rbrace$, thus resulting in $N_d=43$ possible digital states. The number of $\delta$ intervals and the width of each interval are chosen arbitrarily, and there are no restrictions in this respect. The resulting digital states are then sorted to follow the lexicographic order.

\begin{figure}[h!]
\centering
\includegraphics[width=.55\textwidth]{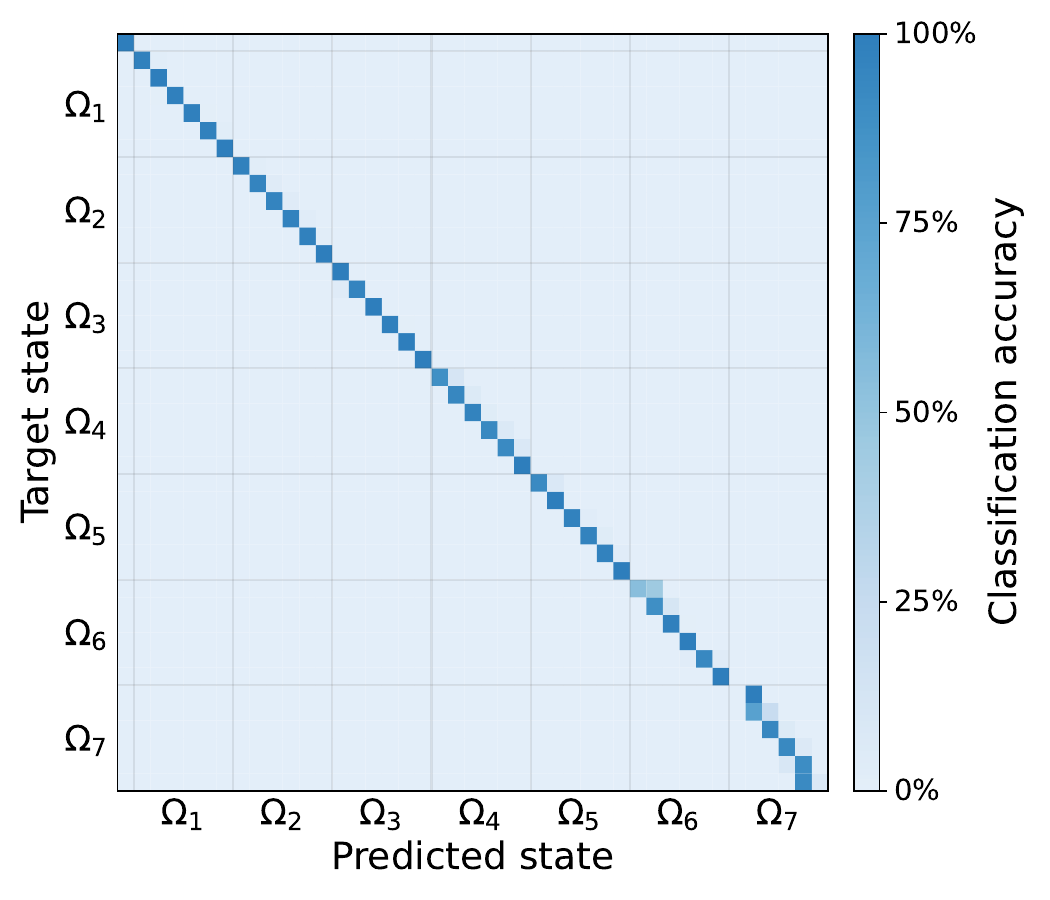}
\caption{L-shaped cantilever beam - Confusion matrix measuring the offline performance of the DL models in correctly categorizing the digital state. Results are reported in terms of classification accuracy, measuring how observational data are classified with respect to the ground truth digital state. Digital states are ordered first for damage location and then for damage level. \label{fig:Lframe_confusion}}
\end{figure}

The confusion matrix reported in \fig\ref{fig:Lframe_confusion} measures the offline performance of $\text{N\hspace{-1px}N}_\text{CL}$ and $\text{N\hspace{-1px}N}_\text{RG}$ in assimilating noisy FOM data to classify the digital state, among the $N_d$ possible outcomes of $D_t$. The (unknown) ground truth digital state is detected by the DL models with an overall classification accuracy of $93.61\%$. Moreover, it can be argued from the confusion matrix that most of the misclassifications are due to the damage scenarios related to a stiffness reduction within $\Omega_6$ or within $\Omega_7$. This is a quite expected outcome since measurements closer to the clamped side are only marginally affected by the presence of damage close to the free end of the beam, thus yielding a smaller sensitivity of sensor recordings to damage. This confusion matrix then serves as the CPT encoding the $\phi_t^\text{N\hspace{-1px}N}$ factor.

For the present case, we consider four possible control inputs, each provided with a CPT modeling the transition probability $p(D_{t+1}|D_t,U^A_t=u^A_t)$ from $D_t$ to $D_{t+1}$ after taking the action $u^A_t$, and collectively encoding the $\phi_t^\text{history}$ factor. These internal models of how structural health
is expected to evolve do not reflect the prescribed ground truth evolution, which is unknown to the DT. The considered control inputs are the following:
\begin{itemize}
\item Do nothing (DN) action. There is no maintenance action planned in this case and the physical state will evolve according to a stochastic deterioration process.
\item Minor imperfect maintenance (MI) action. A maintenance action is performed and the asset may be restored from its current condition to a healthier damage state. This can be traced back to, e.g., patching and sealing cracked surfaces, rectifying and replacing expansion joints, or tightening/replacing loose/missing knot bolts for steel members.
\item Major imperfect maintenance (MA) action. A maintenance action is performed and the asset may be restored from its current condition to a healthier damage state, with a higher probability of improvements than in the previous case. This can be traced back to, e.g., repairing heavily damaged slabs, piers, and steel members, and retrofitting compromised structural elements.
\item Perfect maintenance (PM) action. A maintenance action is performed and the asset is restored from its current condition to the damage-free state. This can be traced back to the replacement of excessively compromised structural elements.
\end{itemize}

\subsubsection{Results: two available actions}

We first illustrate the DT capabilities to assimilate observational data and track the structural health evolution, by restricting the available actions to DN and PM. The (unknown) ground truth evolution of structural health varies depending on the most recently applied control input, which can be either DN or PM. In the absence of maintenance, the physical state evolves following a deterioration process. We prescribe a (simulated) stochastic degradation process that monotonically deteriorates the structural health. The degradation process features a probability of damage inception ($y\neq0$) equal to $0.5$. Damage may develop in any of the predefined regions with $\delta=30\%$, and then propagate with $\delta$ increments sampled from a Gaussian probability density function (pdf) centered at $1.5\%$ and featuring a standard deviation equal to $1\%$ (negative increments are rounded to zero). The effect of a PM action is simulated by restoring the physical state to its undamaged configuration. At each time step during the operation, new observational data are simulated according to the (unknown) ground truth structural health and the most recently enacted control input. The DT assimilates the data and estimates the digital state, eventually suggesting the next control input to enact. Note that the prescribed trajectory of the structural health parameters is arbitrarily chosen to fully display the capabilities of the DT. Nevertheless, the DT would be equally capable of tracking the structural health evolution also considering either more or less aggressive degradation processes.

The state transition model encoding $\phi_t^\text{history}$ is conditioned on the most recently issued control input. The transition probability $p(D_{t+1}|D_t,U^A_t=u^A_t)$ from $D_t$ to $D_{t+1}$ associated with the DN action assumes that damage may start in any subdomain $\Omega_j$, with $j=1,\ldots,N_y$, with probability $0.05$, and then grow to the next $\delta$ interval with the same probability. The transition model assumed for the PM action instead maps the $D_t$ belief to a belief $D_{t+1}$ associated with a damage-free condition, independently of the current condition. The corresponding CPTs are transition matrices, where the diagonal entries represent the probability of staying in the same state. The lower-left and upper-right triangles are associated with the probabilities of the system of deteriorating and improving its condition, respectively. Therefore, the DN transition matrix is a lower-left triangular matrix, with the highest probability assigned to remaining in the same state, consistent with what is expected for the deterioration of civil structures. The transition to the next $\delta$ interval is the second most likely transition, while improvements have a zero probability. Once the structure has reached the last $\delta$ interval, it remains in this condition with a probability equal to $1$. In contrast, the PM transition matrix is an upper-right triangular matrix with probabilities equal to $1$ in the first row.

At each time step, the DT selects a control input $u^A_t\in\lbrace\text{DN, PM}\rbrace$ to be enacted on the asset. Taking a DN action yields a positive reward, but also gives the chance of worsening the asset's structural health. On the other hand, the PM action responds to the degrading structural health, yet yields a negative reward. The computation of the costs associated with the health state and control inputs encapsulates the evaluation of the $\phi_t^\text{reward}$ factor quantifying the performance of the asset. In particular, the two reward functions in \eq\eqref{eq:reward} are defined as:
\begin{equation}
R_t^\text{control}(u^A_t)=\left\{
\begin{array}{ll}
+12,& \text{if $u^A_t=\text{DN}$},\\
-20,& \text{if $u^A_t=\text{PM}$},
\end{array}\right.\; \,\,\,
R_t^\text{health}(d_t)=\left\{
\begin{array}{ll}
+0.1,&\text{if $y=0$},\\
-\text{exp}(6\delta/5),&\text{if $y\in\lbrace1,2,3,4\rbrace$},\\
-\text{exp}(\delta),&\text{if $y\in\lbrace5,6,7\rbrace$},
\end{array}\right.
\end{equation}
where $R_t^\text{control}$ targets the cost assigned to each control input and $R_t^\text{health}$ measures the cost associated with the structural health state. These non-dimensional rewards represent indicative values the decision-maker is charged due to the condition of the structure. Although these values are not based on real data, actual values are not usually hard to find. State agencies and companies provide lists with services and costs \cite{art:papa_shinozuka}. The three cases in $R_t^\text{health}$ distinguish between the absence of damage, the presence of damage within the harm closed to the clamped side, and the presence of damage far from the clamp, respectively. Note how these penalize the progressive deterioration of the structural health as a function of $\delta$. $R_t^\text{health}$ can resemble a variety of aspects, like reduction in the level of service due to deterioration, working accidents, structural reliability, and structural failure probability~\cite{art:papa_shinozuka}.

During the offline phase, we solve the planning problem induced by the PGM to compute the control policy $\pi(D_t)$, which maps the digital state belief to actions and encodes the control factor $\phi_t^\text{control}$. The optimization of $R_t(U_t,D_t)$ is carried out as described in \sez\ref{sec:pgm}, assuming a discount factor $\gamma=0.95$ and a weighting factor $\alpha=2$. The computed control policy $\pi(D_t)$ recommends that the asset operates until when $\delta\in[65\%,75\%]$ and $\delta\in[75\%,80\%]$, respectively if $y\in\lbrace1,2,3,4\rbrace$ and if $y\in\lbrace5,6,7\rbrace$, at which point it should be repaired.

\begin{figure}[h!]
\centering
\includegraphics[width=.75\textwidth, trim=30 10 110 20, clip]{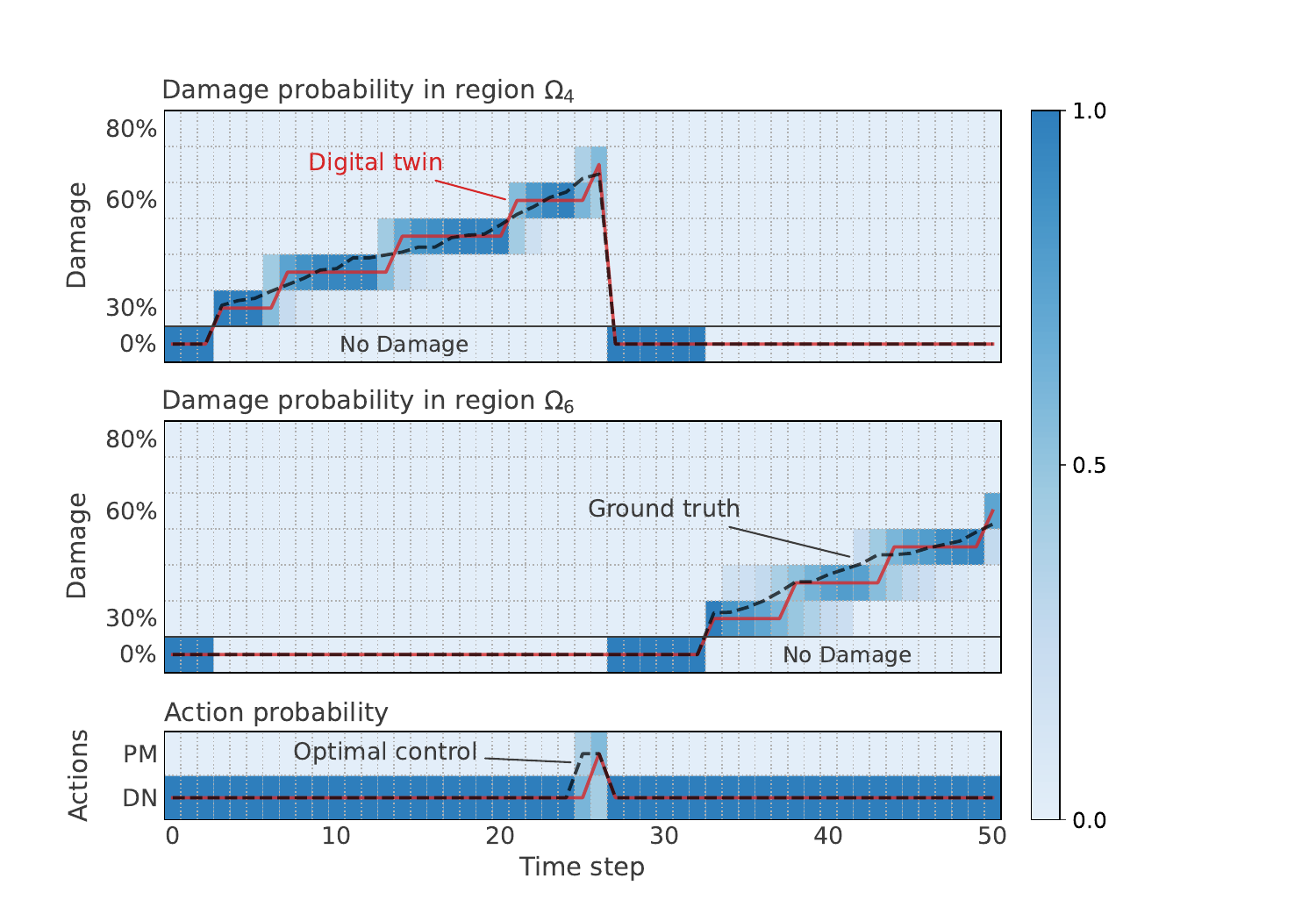}
\caption{L-shaped cantilever beam - Online phase of the digital twin framework with two possible actions: DN (do nothing), and PM (perfect maintenance). Probabilistic and best point estimates of: (top) digital state evolution against the ground truth digital state; (bottom) control inputs informed by the digital twin, against the optimal control input under ground truth. In the top panels the background color corresponds to $p(D_t | D_{t-1}, D^{\text{N\hspace{-1px}N}}_t, U^A_{t-1}=u^A_{t-1})$. In the bottom panel it corresponds to $p(U_t | D_t)$.\label{fig:history_2act}}
\end{figure}

\fig\ref{fig:history_2act} depicts a simulated online phase of the DT up to time step $t_c=50$. Results are reported in terms of the (unknown) ground truth digital state, and the corresponding DT estimate after assimilating the observational data. The graphs report the evolution of the digital state only for the damaged regions, nevertheless, damage can potentially affect all $\Omega_1,\ldots,\Omega_7$ predefined damaged regions. The DT proves capable of accurately tracking the digital state evolution with relatively low uncertainty. The corresponding estimation of the control inputs is reported in the bottom part of the figure, demonstrating that the DT is able to promptly suggest the PM action within one time step of when the (unknown) ground truth structural health demands it.

\begin{figure}[h!]
\centering
\includegraphics[width=.75\textwidth, trim=30 0 110 20, clip]{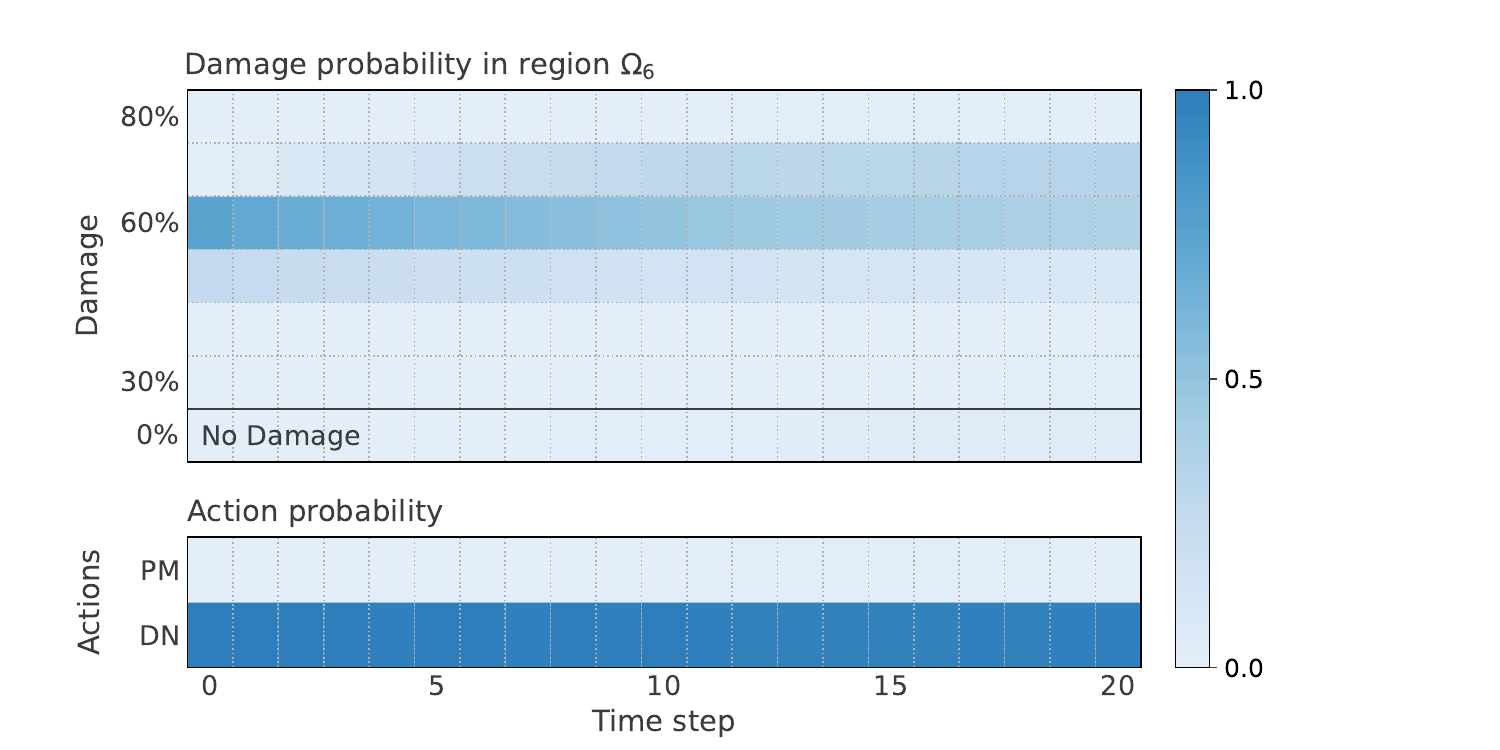} 
\caption{L-shaped cantilever beam - Digital twin future predictions with two possible actions: DN (do nothing), and PM (perfect maintenance). The starting time is $t_c=50$. In the top panel the probability $p(D_t | D_{t-1}, U_{t-1})$ relates to the amount of damage in $\Omega_6$. In the bottom panel it corresponds to $p(U_t | D_t)$.\label{fig:predictions_2act}}
\end{figure}

\fig\ref{fig:predictions_2act} depicts the predicted evolution of the digital state and of the corresponding informed control inputs, starting from $t_c=50$. The prediction horizon is extended over $20$ time steps in the future so that $t_p=t_c+20$. The DT prediction engine informs about the expected future degradation of the structural health, allowing to plan future interventions.

\subsubsection{Results: four available actions}
We now consider all four possible control inputs. We prescribe a stochastic degradation process with a probability of damage inception ($y\neq0$) equal to $0.5$. Damage may develop in any of the predefined regions with damage level sampled from a uniform distribution $\delta\in[30\%,70\%]$, and then propagate as in the previous case. This more aggressive degradation process is used to spot in a few time steps the effectiveness of the decision-making capabilities of the DT.
The effect of the MI and MA actions on the asset is simulated according with stochastic repair processes, for which the structural health is forced to monotonically improve. The effect of a MI action is simulated with $\delta$ decrements sampled from a Gaussian pdf centered at $-12.5\%$ and featuring a standard deviation equal to $1\%$, while the effect of a MA action is modeled with $\delta$ decrements sampled from a Gaussian pdf centered at $-17.5\%$ and featuring a standard deviation equal to $1\%$. In both cases, the damage-free condition is assumed to be recovered if the resulting structural state features $\delta<30\%$.

The transition model $p(D_{t+1}|D_t,U^A_t=u^A_t)$ associated with the MI action assumes no improvement in the structural health with probability $0.1$, improvement of one $\delta$ interval with probability $0.75$, and improvement of two $\delta$ intervals with probability $0.15$. The resulting CPT is an upper-right triangular transition matrix, as deterioration from any state upon a repair action is assumed to have zero probability. The highest probability is assigned to improvements of one $\delta$ interval, followed by improvements of two $\delta$ intervals. There is also a lower probability of remaining in the same deteriorated state, which reflects a failed maintenance. Similarly, the MA action assumes no improvement with probability $0.05$, improvement of one $\delta$ interval with probability $0.3$, improvement of two $\delta$ intervals with probability $0.4$, and improvement of three $\delta$ intervals with probability $0.25$. In this case, the highest probability is assigned to improvements of two $\delta$ intervals, followed by improvements of one $\delta$ intervals, three $\delta$ intervals, and finally, the lowest probability is associated with the possibility of a failed maintenance.

The two reward functions in \eq\eqref{eq:reward} are chosen as:
\begin{equation}
R_t^\text{control}(u^A_t)=\left\{
\begin{array}{ll}
+12,& \text{if $u^A_t=\text{DN}$},\\
-20,& \text{if $u^A_t=\text{PM}$},\\
-8,& \text{if $u^A_t=\text{MI}$},\\
-15,& \text{if $u^A_t=\text{MA}$},
\end{array}\right.\; \,\,\,
R_t^\text{health}(d_t)=\left\{
\begin{array}{ll}
+0.1,&\text{if $y=0$},\\
-\text{exp}(5\delta),&\text{if $y\in\lbrace1,2,3,4\rbrace$},\\
-\text{exp}(4\delta),&\text{if $y\in\lbrace5,6,7\rbrace$}.
\end{array}\right.
\end{equation}
We assume a discount factor $\gamma=0.95$, and a weighting factor $\alpha=2.5$. The resulting control policy $\pi(D_t)$ recommends that the asset should operate until when $\delta\in[30\%,35\%]$, after which: if $y\in\lbrace1,2,3,4\rbrace$, a MI action should be performed when $\delta\in[35\%,45\%]$, and a PM action should be performed when $\delta>45\%$; while, if $y\in\lbrace5,6,7\rbrace$, the MI and MA actions should be performed, respectively when $\delta\in[35\%,55\%]$ and when $\delta\in[55\%,75\%]$, and a PM action should be performed when $\delta>75\%$.

\begin{figure}[h!]
\centering
\includegraphics[width=.75\textwidth, trim=30 40 110 20, clip]{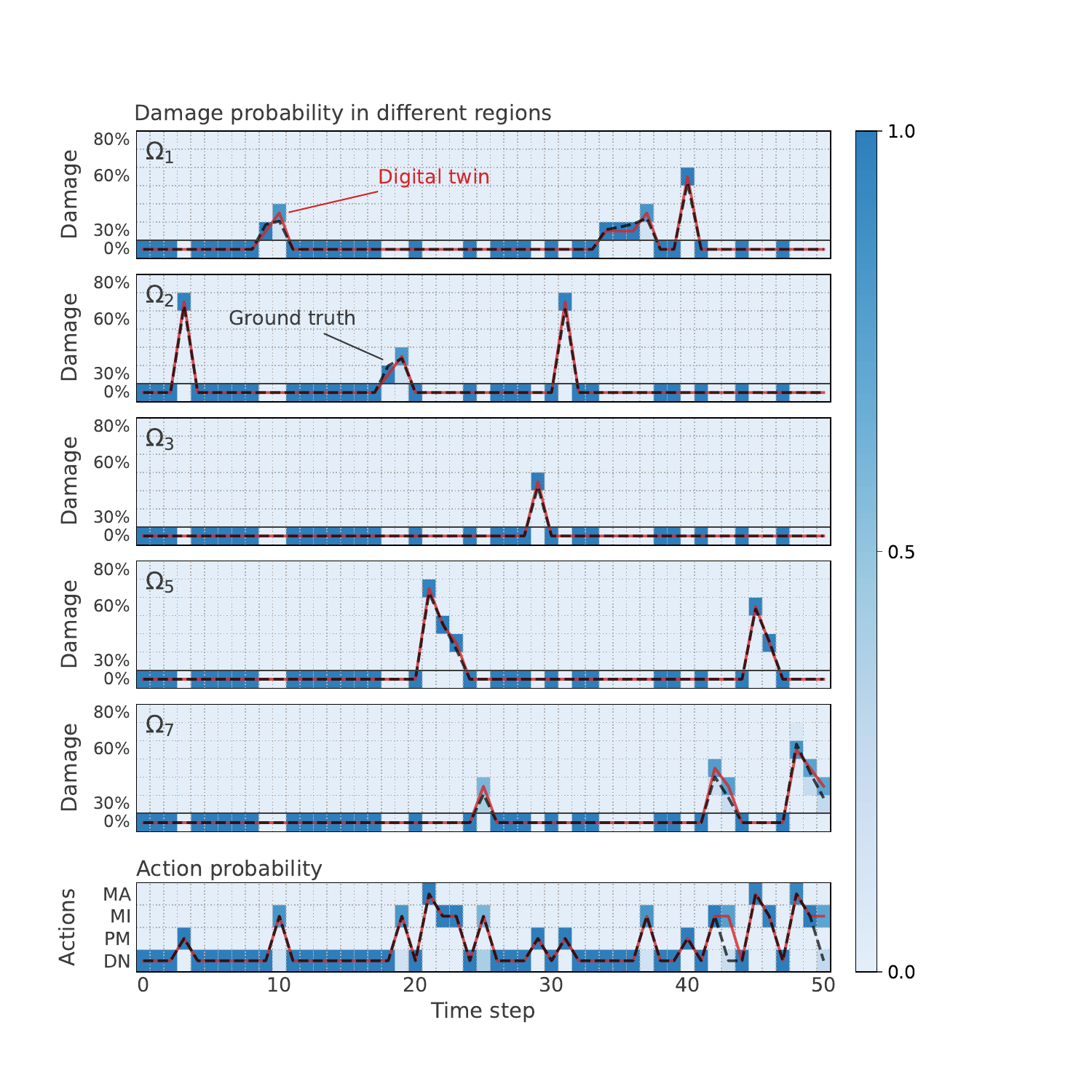}
\caption{L-shaped cantilever beam - Online phase of the digital twin framework with four possible actions: DN (do nothing), PM (perfect maintenance), MI (minor imperfect maintenance), and MA (major imperfect maintenance). Probabilistic and best point estimates of: (top) digital state evolution against the ground truth digital state; (bottom) control inputs informed by the digital twin, against the optimal control input under ground truth. In the top panels the background color corresponds to $p(D_t | D_{t-1}, D^{\text{N\hspace{-1px}N}}_t, U^A_{t-1}=u^A_{t-1})$. In the bottom panel it corresponds to $p(U_t | D_t)$.\label{fig:history_4act}}
\end{figure}

\fig\ref{fig:history_4act} depicts a simulated online phase of the DT up to $t_c=50$. The DT accurately tracks the digital state evolution and timely suggests the appropriate control inputs most of the time. In particular, the DT proposes the optimal control input, except for the time steps $t=43$ and $t=50$ featuring a sub-optimal action. In both cases, a MI action is proposed in place of a DN, because the DT estimates a $\delta\in[35\%,45\%]$ instead of a $\delta\in[30\%,35\%]$ related to a stiffness reduction within $\Omega_7$. This is in line with what was observed in the confusion matrix of \fig\ref{fig:Lframe_confusion}, due to the limited sensitivity of recordings to damage scenarios affecting the terminal region of the beam. This peculiar type of misclassification turns out to be the most pathological in the confusion matrix and is therefore capable of potentially spoiling the assimilation of observational data. Nevertheless, the DT reverts to correctly tracking the structural health of the asset within one time step.

\begin{figure}[h!]
\begin{centering}
\includegraphics[width=.75\textwidth, trim=30 0 110 20, clip]{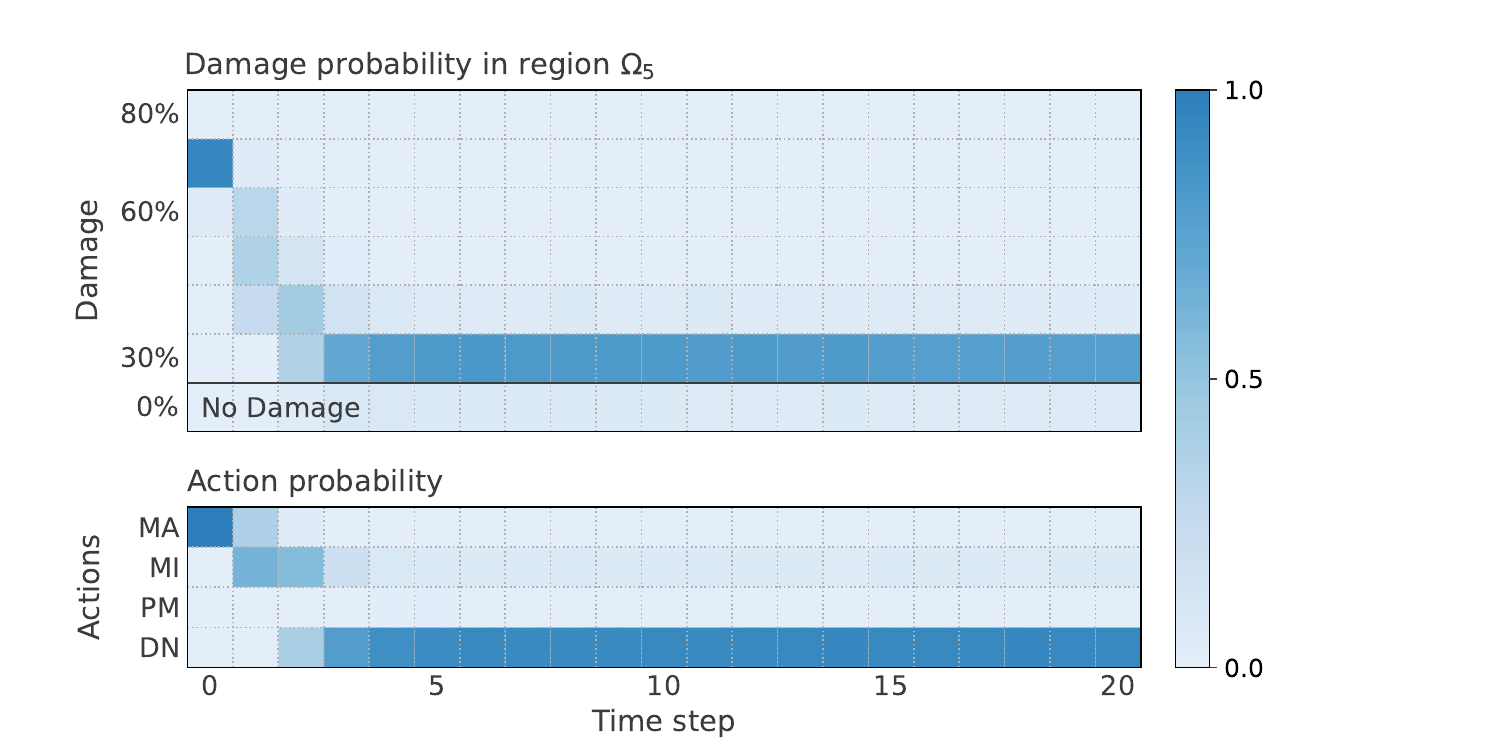}
\caption{L-shaped cantilever beam - Digital twin future predictions with four possible actions: DN (do nothing), PM (perfect maintenance), MI (minor imperfect maintenance), and MA (major imperfect maintenance). The starting time is $t_c=21$. In the top panel the probability $p(D_t | D_{t-1}, U_{t-1})$ relates to the amount of damage in $\Omega_5$. In the bottom panel it corresponds to $p(U_t | D_t)$.\label{fig:predictions_4act}}
\end{centering}
\end{figure}

\fig\ref{fig:predictions_4act} depicts the predicted evolution of the digital state and control inputs, from $t_c=21$ and over $20$ time steps in the future. The DT prediction correctly suggests taking with high probability a MA action, followed by two MI actions, and accordingly predicts the corresponding evolution of the structural health. Comparing the DT prediction with what is effectively experienced during the online phase (see \fig\ref{fig:history_4act}), note how the DT prediction closely resembles the actual evolution of the digital state and control inputs. This is a remarkable result in terms of DT prediction capabilities, since the DT is not aware of the future values of the structural health parameters, and the relative transition models do not match their real (stochastic) evolution.

\subsection{Railway bridge}
\label{sec:bridge}

\begin{figure}[h!]
\begin{center}
\includegraphics[width=.8\textwidth]{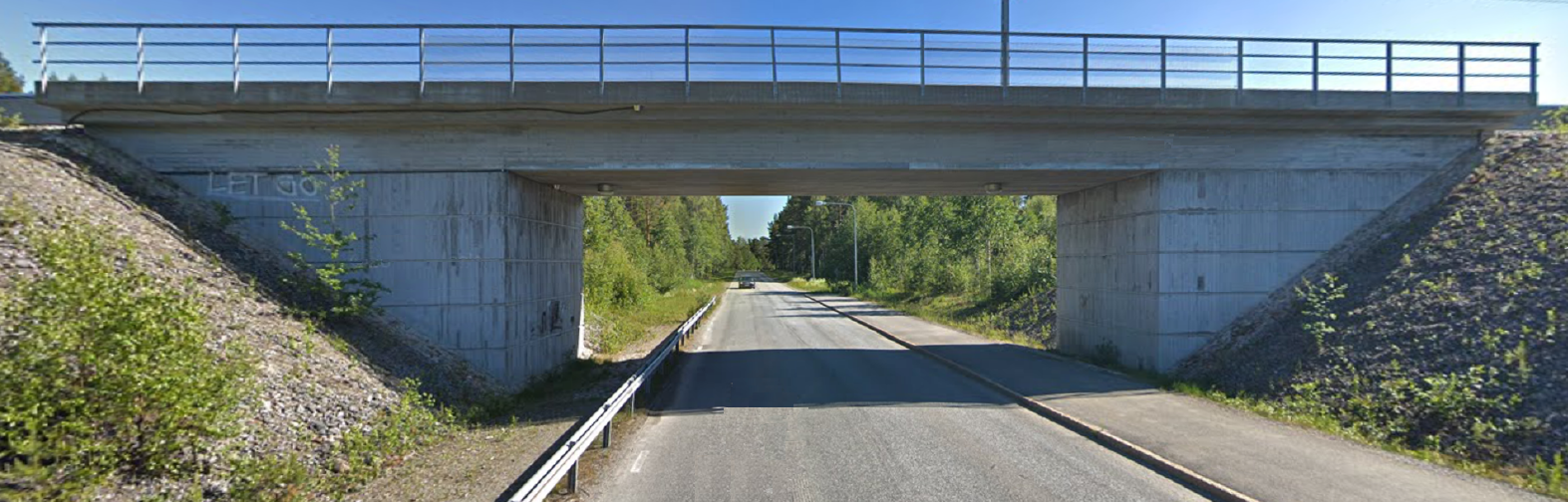}
\caption{H{\"o}rnefors railway bridge.\label{fig:bridge_photo}}
\end{center}
\end{figure}

The second case study concerns the railway bridge depicted in \fig\ref{fig:bridge_photo}. It is an integral concrete portal frame bridge located along the Bothnia line in the Swedish suburbs of H{\"o}rnefors. It features a span of $15.7~\textup{m}$, a free height of $4.7~\textup{m}$, and a width of $5.9~\textup{m}$ (edge beams excluded). The thickness of the structural elements is $0.5~\textup{m}$ for the deck, $0.7~\textup{m}$ for the frame walls, and $0.8~\textup{m}$ for the wing walls. The bridge is founded on two plates connected by stay beams and supported by pile groups. The concrete is of class C35/45, whose mechanical properties are: $E=34~\textup{GPa}$, $\nu= 0.2$, $\rho=2500~\textup{kg/m}^3$. The superstructure consists of a single track with sleepers spaced $0.65~\textup{m}$ apart, resting on a ballast layer $0.6~\textup{m}$ deep, $4.3~\textup{m}$ wide and featuring a density $\rho_{B}=1800~\textup{kg/m}^3$. The geometrical and mechanical modeling data have been adapted from former research activities on the relevant soil-structure interaction, see~\cite{thesis:kth3,thesis:kth2}.

The bridge is subjected to the transit of Gr{\"o}na T{\r a}get trains type, at a speed $\upsilon\in[160,215]~\textup{km/h}$. Only trains composed of two wagons are considered, thus characterized by $8$ axles, each one carrying a mass $\psi\in[16,22]~\textup{ton}$. The corresponding load model is described in~\cite{art:Metodologico}, and consists of $25$ equivalent distributed forces transmitted by the sleepers to the deck through the ballast layer with a slope $4:1$, according with Eurocode 1~\cite{code:EC1}.

\subsubsection{Dataset assembly}

\begin{figure}[h!]
\begin{centering}
\includegraphics[width=.6\textwidth]{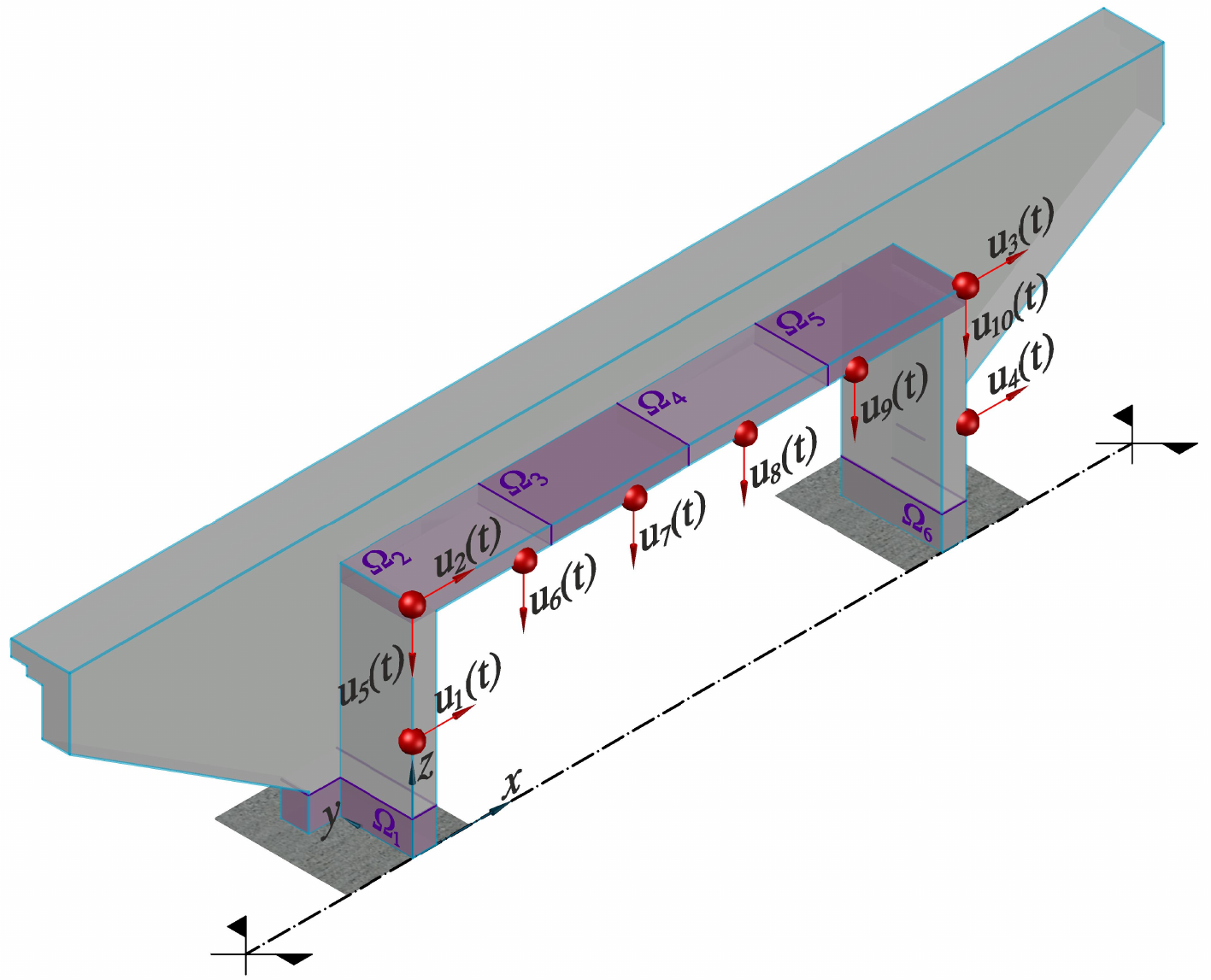}
\caption{Railway bridge: details of synthetic recordings related to displacements $u_1(t),\ldots,u_{10}(t)$, and predefined damage regions $\Omega_1,\ldots,\Omega_6$.\label{fig:Ponte_Model}}
\end{centering}
\end{figure}

Synthetic displacement time histories $\mathbf{U}$ are obtained from $N_u=10$ sensors deployed as depicted in \fig\ref{fig:Ponte_Model}. Each recording is provided for a time interval $(0,T=1.5~\textup{s})$ with an acquisition frequency $f_\text{s}=400~\textup{Hz}$. This setting allows to record train passages at the lowest speed of $160~\textup{km/h}$, and properly catches the structural response at the maximum speed of $215~\textup{km/h}$. Recordings are corrupted with an additive Gaussian noise yielding a signal-to-noise ratio of $120$. 

In addition to the undamaged condition, the presence of damage in the structure is accounted for using a localized stiffness reduction that can take place within $N_y=6$ predefined subdomains $\Omega_j$, with $j=1,\ldots,N_y$, as depicted in \fig\ref{fig:Ponte_Model}. The stiffness reduction can occur with a magnitude $\delta\in[30\%,80\%]$, and is kept fixed while a train travels across the bridge. 

The FOM features $N_\text{FE}=17,292$ dofs, resulting from a finite element discretization with an element size of $0.80~\textup{m}$ and a reduced size of $0.15~\textup{m}$ for the deck, to enable a smooth propagation of the traveling load. The presence of the ballast layer is accounted for through an increased density for the deck and for the edge beams. The embankments are accounted for through distributed springs, modeled as a Robin mixed boundary condition (with elastic coefficient $k_{\textup{robin}}=10^{8}~\textup{N/m}^3$) applied on the surfaces facing the ground. The structural dissipation is modeled by means of a Rayleigh's damping matrix, assembled to account for a $5\%$ damping ratio on the first two structural modes. 

The ROM is obtained from a snapshot matrix $\mathbf{S}$, assembled through $400$ evaluations of the FOM for different values of parameters $\boldsymbol{\mu} = (\upsilon,\psi,y,\delta)^\top$. By setting the error tolerance to $\epsilon=10^{-3}$, $N_\text{RB}=133$ POD modes are to be considered.

The training dataset $\mathcal{D}$ is built with $I=10,000$ instances collected using the ROM. Also in this case, the testing phase of $\text{N\hspace{-1px}N}_\text{CL}$ and of $\text{N\hspace{-1px}N}_\text{RG}$ is carried out considering noisy FOM solutions. The monitoring of the asset is then simulated by assimilating $N_\text{obs}=1$ noisy observations at each time step. As the structural health of the bridge evolves over time, the DT estimates the variation in the structural health parameters every time a train travels across the bridge.

\subsubsection{Digital twin framework}

As in the previous case, the two structural health parameters within the digital state are $\boldsymbol{d}=(y,\delta)^\top$. The range in which the damage level $\delta$ can take values is discretized in $N_\delta=6$ intervals. The resulting $N_d=37$ possible digital states are sorted first for damage location and then for damage level.

\begin{figure}[h!]
\begin{centering}
\includegraphics[width=.55\textwidth]{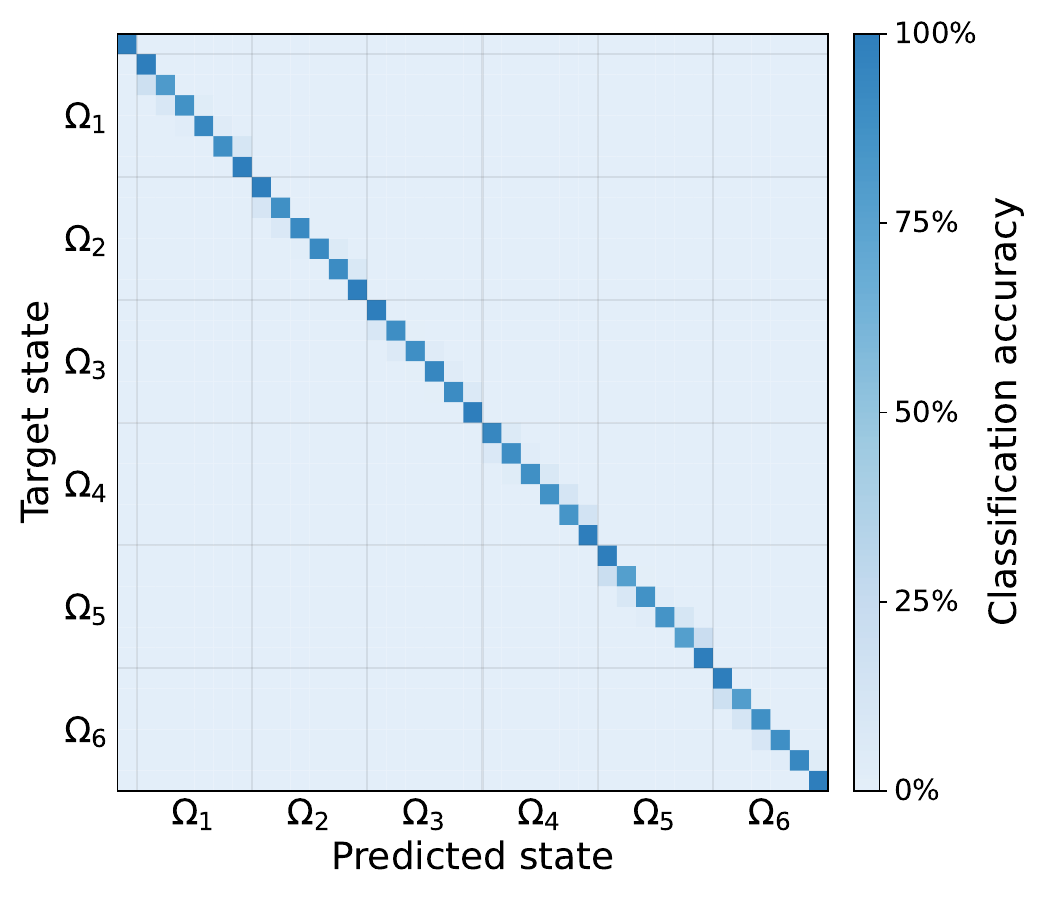}
\caption{Railway bridge - Confusion matrix measuring the offline performance of the DL models in correctly categorizing the digital state. Results are reported in terms of classification accuracy, measuring how observational data are classified with respect to the ground truth digital state. Digital states are ordered first for damage location and then for damage level. \label{fig:brige_confusion}}
\end{centering}
\end{figure}

The confusion matrix measuring the offline performance of $\text{N\hspace{-1px}N}_\text{CL}$ and of $\text{N\hspace{-1px}N}_\text{RG}$ in correctly categorizing the digital state is reported in \fig\ref{fig:brige_confusion}. The ground truth digital state is detected with an overall classification accuracy of $91.39\%$. In this case, the majority of misclassifications are due to confusing adjacent digital states relative to the same damage location, thus yielding a tridiagonal band matrix.

For the present case, we consider the following three possible control inputs:
\begin{itemize}
\item Do nothing (DN) action. There is no maintenance action planned in this case and the physical state will evolve according to a stochastic deterioration process.
\item Perfect maintenance (PM) action. A maintenance action is performed and the asset is restored from its current condition to the damage-free state.
\item Restrict operational conditions (RE) action. The operational conditions of the bridge are restricted by allowing only lightweight trains, carrying less than $18~\textup{ton}$ per axle, to travel across the bridge. Such a restriction results in a lower deterioration rate, but also yields a lower revenue generated by the infrastructure.
\end{itemize}

In the cases where the most recently issued control input is either DN or RE, the physical state undergoes a degradation process that monotonically deteriorates the structural health. When operational conditions are not restricted, we prescribe a stochastic degradation process featuring a probability of damage inception ($y\neq0$) equal to $0.5$. Damage may develop in any of the predefined regions with damage level sampled from a uniform distribution $\delta\in[30\%,35\%]$, and then propagate with $\delta$ increments sampled from a Gaussian pdf centered at $1.5\%$ and featuring a standard deviation equal to $1\%$ (negative increments are rounded to zero). When the operations are restricted and only lightweight trains are allowed to travel across the bridge, we instead assume a probability of damage inception equal to $0.25$. In this eventuality, damage may develop with damage level sampled from a uniform distribution $\delta\in[30\%,35\%]$, and then propagate with $\delta$ increments sampled from a Gaussian pdf centered at $0.95\%$ and featuring a standard deviation equal to $0.5\%$. The resulting trajectory of the structural health parameters is intended to represent periods of gradual degradation in the structural health, as well as sudden changes due to discrete damage events. Also in this case, the effect of a PM action is simulated by restoring the physical state to its undamaged configuration.

The transition model $p(D_{t+1}|D_t,U^A_t=u^A_t)$ associated with the DN action assumes that damage may start in any subdomain $\Omega_j$, with $j=1,\ldots,N_y$, with probability $0.1$, and then grow to the next $\delta$ interval with the same probability. For the transition model associated with the RE action, this probability is assumed to decrease to $0.03$. The CPTs associated with the DN and RE actions are therefore lower-left triangular transition matrices. The highest probability assigned to remaining in the same state, followed by the transition to the next $\delta$ interval, with zero probability of improvements. The transition model assumed for the PM action instead maps the $D_t$ belief to a belief $D_{t+1}$ associated with a damage-free condition, independently of the current condition. The CPT associated with the PM action is therefore an upper-right triangular transition matrix with probabilities equal to $1$ in the first row.
            
In this case, the two reward functions in \eq\eqref{eq:reward} are chosen as:
\begin{equation}
R_t^\text{control}(u^A_t)=\left\{
\begin{array}{ll}
+30,& \text{if $u^A_t=\text{DN}$},\\
-250,& \text{if $u^A_t=\text{PM}$},\\
+27,& \text{if $u^A_t=\text{RE}$},
\end{array}\right.\; \,\,\,
R_t^\text{health}(d_t)=\left\{
\begin{array}{ll}
+0,&\text{if $y=0$},\\
-\text{exp}(5\delta)+4,&\text{if $y\neq0$},\\
-250,&\text{if $\delta\geq79\%$},
\end{array}\right.
\end{equation}
where the last contribution in $R_t^\text{health}$ penalizes excessively compromised structural states with a significantly negative reward.

\subsubsection{Results}

During the offline phase, we solve the planning problem in \eq\eqref{eq:optimization} by assuming a discount factor $\gamma=0.90$, and a weighting factor $\alpha=1$. The resulting control policy $\pi(D_t)$ recommends that the asset operates in ordinary conditions until when $\delta\in[30\%,
35\%]$, after which point it should fall back to the more conservative RE regime in order to minimize further degradation. Once reached $\delta\geq65\%$, the bridge should be finally repaired.

\begin{figure}[h!]
\begin{centering}
\includegraphics[width=.75\textwidth, trim=30 10 110 20, clip]{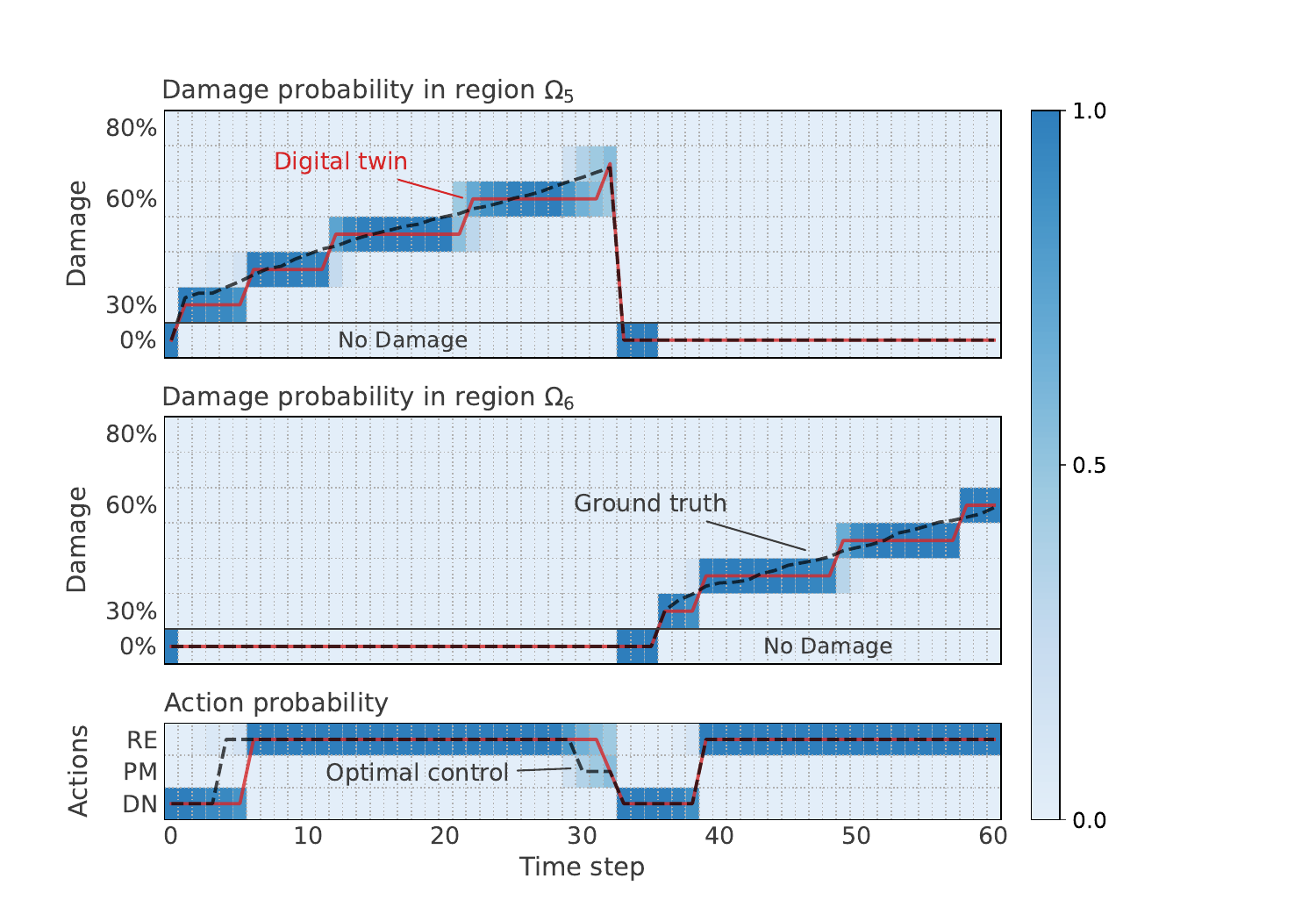}
\caption{Railway bridge - Online phase of the digital twin framework with three possible actions: DN (do nothing), PM (perfect maintenance), and RE (restrict operational conditions). Probabilistic and best point estimates of: (top) digital state evolution against the ground truth digital state; (bottom) control inputs informed by the digital twin, against the optimal control input under ground truth. In the top panels the background color corresponds to $p(D_t | D_{t-1}, D^{\text{N\hspace{-1px}N}}_t, U^A_{t-1}=u^A_{t-1})$. In the bottom panel it corresponds to $p(U_t | D_t)$.\label{fig:history_3act}}
\end{centering}
\end{figure}

\fig\ref{fig:history_3act} reports a sample simulation of the DT online phase up to time step $t_c=60$. The DT correctly tracks the digital state with relatively low uncertainty. Damage initially develops within $\Omega_5$, and the DT follows its evolution with a limited delay of at most two time steps, with respect to the ground truth, due to the need of updating the relative prior belief from the previous time steps. The RE action is suggested as soon as the DT estimates a $\delta\in[35\%,65\%]$, after which point the DT keeps on tracking the structural health parameters evolving with a lower deterioration rate. A PM action is finally suggested due to an excessively compromised structural state. A similar behavior can be observed for the following damage scenario affecting $\Omega_6$.

\begin{figure}[h!]
\centering
\includegraphics[width=.75\textwidth, trim=30 0 110 20, clip]{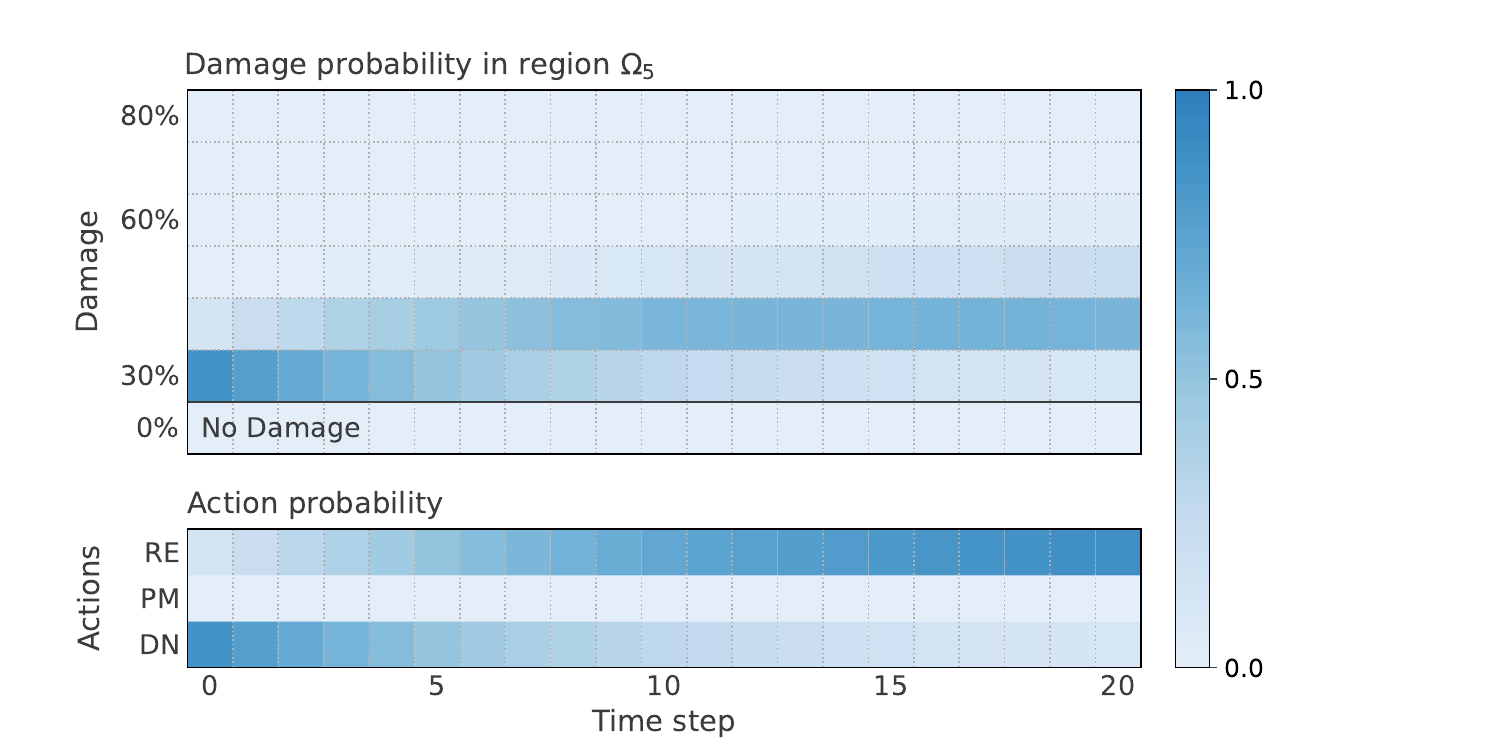}
\caption{Railway bridge - Digital twin future predictions with three possible actions: DN (do nothing), PM (perfect maintenance), and RE (restrict operational conditions). The starting time is $t_c=5$. In the top panel the probability $p(D_t | D_{t-1}, U_{t-1})$ relates to the amount of damage in $\Omega_5$. In the bottom panel it corresponds to $p(U_t | D_t)$.\label{fig:predictions_3act}}
\end{figure}

\fig\ref{fig:predictions_3act} reports the predicted evolution of the digital state and control inputs, from $t_c=5$ and over $20$ time steps in the future. The DT predicts the expected degradation of the structural health according to the transition model associated with the DN action, before predicting to take a RE action with relatively high probability after a few time steps. The DT prediction is close to what is effectively experienced online (see \fig\ref{fig:history_3act}). However, besides having the estimated digital state two time steps behind the ground truth value, the prediction is also too optimistic in terms of deterioration rate, which suggests the use of a more refined transition model.

\section{Conclusions}
\label{sec:conclusions}
In this work we have proposed a predictive digital twin approach to the health monitoring, maintenance, and management planning of civil structures, to advance condition-based and predictive maintenance practices. The presented strategy relies upon a probabilistic graphical model inspired by~\cite{art:kapteyn2021probabilistic}. This framework is used to encode the asset-twin coupled dynamical system, the relevant end-to-end information flow via observational data (physical to digital) and control inputs (digital to physical), and its evolution over time, all with quantified uncertainty. The assimilation of observational data is carried out with deep learning models, leveraging the capabilities of convolutional layers to  automatically select and extract damage-sensitive features from raw vibration recordings. The structural health parameters comprising the digital state are used to capture the variability of the physical asset. They are continually updated in a sequential Bayesian inference fashion, according to control-dependent transition dynamics models describing how the structural health is expected to evolve. The updated digital state is eventually exploited to predict the future evolution of the physical system and the associated uncertainty. This enables predictive decision-making about maintenance and management actions. 

The computational procedure takes advantage of a preliminary offline phase which involves: (i) using physics-based numerical models and reduced order modeling, to overcome the lack of experimental data for civil applications under varying damage and operational conditions while populating the datasets for training the deep learning models; (ii) learning the health-dependent control policy to be applied at each time step of the online phase, to map the belief over the digital state onto actions feeding back to the physical asset.

The proposed strategy has been assessed against the simulated monitoring of an L-shaped cantilever beam and a railway bridge. In the absence of experimental data, the tests have been carried out considering high-fidelity simulation data, corrupted with an additive Gaussian noise. The obtained results have proved the digital twin capabilities of accurately tracking the digital state evolution under varying operational conditions, with relatively low uncertainty. The framework is also able to promptly suggest the appropriate control input, within at most two time steps of when the (unknown) ground truth structural health demands it.

Although the capabilities of health-aware digital twins are showcased in the specific context of monitoring the structural integrity of civil structures to advance predictive maintenance practices, the applicability of the presented framework is general. Indeed, the proposed framework can be adapted for various types of structures and engineering systems by adjusting the components within the dynamic Bayesian network to align with the specific characteristics of the problem at hand. The solution to the inverse problem (if any) can be estimated by assimilating available observational data using methods other than deep neural networks, for instance through Markov chain Monte Carlo sampling algorithms. Similarly, the state transition models are closely tied to the employed parametrization of the digital state and the availability of historical data. The same applies to the available control inputs, which are likely to vary for different structures, such as those in mechanical or aerospace systems, and the method chosen for solving the associated planning problem. Additionally, the graph topology can be easily reorganized to adapt to situations where observational data are not acquired after issuing a control input, or when control inputs are issued with a different frequency than that governing the digital twin update.

Future research lines will investigate the ability of the digital twin to update the transition dynamics models by learning from previous data. As suggested by the railway bridge case study, this will allow for a more accurate prediction of the expected evolution of the digital state, thus enabling predictive decision-making better tailored to the monitored asset. Another aspect of interest concerns solving the planning problem induced by the probabilistic model using reinforcement learning algorithms, capable of taking into account a finite planning horizon representing the design lifetime of the asset.

\section*{Data availability}
The observational data used to run the experiments presented in \sez\ref{sec:experiments} are available in the public repository \texttt{digital-twin-SHM}~\cite{Repo}. The \texttt{Matlab} library for finite element simulation and reduced-order modeling of partial differential equations employed to generate these data is available in the public repository \texttt{Redbkit}~\cite{Redbkit}.

\section*{Code availability}
The implementation code used for the experiments presented in \sez\ref{sec:experiments} is available in the public repository \texttt{digital-twin-SHM}~\cite{Repo}. The code implements the proposed digital twin framework and can be used to generate the graphs of digital state estimation and prediction reported in this paper. The DL models trained according to the implementation details reported in Appendix~\ref{sec:implementation} are also made available in the same repository.

\section*{Acknowledgements}
Matteo Torzoni gratefully acknowledges the financial support from the Politecnico di Milano through the interdisciplinary Ph.D. Grant ‘‘Physics-informed deep learning for structural health monitoring”. Marco Tezzele and Karen E. Willcox acknowledge support from the NASA University Leadership Initiative under Cooperative Agreement 80NSSC21M0071. Andrea Manzoni acknowledges the project ‘‘Dipartimento di Eccellenza” 2023-2027, funded by MUR, and the project FAIR (Future Artificial Intelligence Research), funded by the NextGenerationEU program within the PNRR-PE-AI scheme (M4C2, Investment 1.3, Line on Artificial Intelligence).

\bibliographystyle{ieeetr}

\appendix
\newpage
\section{Implementation details}
\label{sec:implementation}

In this Appendix, we discuss the implementation details of the DL models described in \sez\ref{sec:DL_models}. The architectures, as well as the relevant hyperparameters and training options, have been chosen through a preliminary study, aimed at minimizing $\mathcal{L}_\text{CL}$ and $\mathcal{L}_\text{RG}$, while retaining the generalization capabilities of $\text{N\hspace{-1px}N}_\text{CL}$ and of $\text{N\hspace{-1px}N}^j_\text{RG}$, with $j=1,\ldots,N_y$. Since all $\text{N\hspace{-1px}N}^j_\text{RG}$ models share the same architecture, the index $j$ will be dropped in the following for ease of notation.

In the present work, $\text{N\hspace{-1px}N}_\text{CL}$ and $\text{N\hspace{-1px}N}_\text{RG}$ are set as 12-layers DL models, whose architecture is outlined in \tab\ref{tab:NN_CL_arch} and in \tab\ref{tab:NN_RG_arch}, respectively. $\text{N\hspace{-1px}N}_\text{CL}$ and $\text{N\hspace{-1px}N}_\text{RG}$ feature a damage-sensitive feature extractor required to be insensitive to transformations in the input not related to damage. This is implemented through the composition of three one-dimensional (1D) convolutional units. Convolutional layers naturally embed good relational inductive biases such as locality and translation equivariance~\cite{art:Battaglia,book:DL_book}, which prove highly effective to detect time correlations within and across time series. The resulting sparse connectivity and parameter sharing also make them computationally efficient. Each convolutional unit consists of a convolutional layer, followed by a Tanh activation function, max pooling, and dropout. The extracted features are expected to be sensitive to the presence of damage, but insensitive to measurement noise and operational variability. The extracted features are then reshaped through a flatten layer and run through a stack of three fully-connected layers: the first two are Tanh-activated, while the output layer of $\text{N\hspace{-1px}N}_\text{CL}$ is Softmax-activated, and the output layer of $\text{N\hspace{-1px}N}_\text{RG}$ has no activation function.

\begin{table}[h!]
\caption{$\text{N\hspace{-1px}N}_\text{CL}$ - (a) employed architecture, and (b) selected hyperparameters and training options.}
      \centering
       \subfloat[\label{tab:NN_CL_arch}]{
       \scriptsize
\begin{tabular}{llll}
    \toprule
    \mbox{Layer} & \mbox{Output shape} &  \mbox{Activation} & \mbox{Input}\\
    \toprule
    \mbox{0 - Input} & \mbox{$(B_\text{CL},L,N_u)$} & \mbox{None} & \mbox{--}\\
    \mbox{1 - Conv1D} & \mbox{$(B_\text{CL},L,32)$} & \mbox{Tanh} & \mbox{0}\\
    \mbox{2 - MaxPooling1D} & \mbox{$(B_\text{CL},L/2,32)$} & \mbox{None} & \mbox{1}\\
    \mbox{3 - Dropout} & \mbox{$(B_\text{CL},L/2,32)$} & \mbox{None} & \mbox{2}\\
    \mbox{4 - Conv1D} & \mbox{$(B_\text{CL},L/2,64)$} & \mbox{Tanh} & \mbox{3}\\
    \mbox{5 - MaxPooling1D} & \mbox{$(B_\text{CL},L/4,64)$} & \mbox{None} & \mbox{4}\\
    \mbox{6 - Dropout} & \mbox{$(B_\text{CL},L/4,64)$} & \mbox{None} & \mbox{5}\\
    \mbox{7 - Conv1D} & \mbox{$(B_\text{CL},L/4,32)$} & \mbox{Tanh} & \mbox{6}\\
    \mbox{8 - MaxPooling1D} & \mbox{$(B_\text{CL},L/8,32)$} & \mbox{None} & \mbox{7}\\
    \mbox{9 - Dropout} & \mbox{$(B_\text{CL},L/8,32)$} & \mbox{None} & \mbox{8}\\
    \mbox{10 - Flatten} & \mbox{$(B_\text{CL},4L)$} & \mbox{None} & \mbox{9}\\
    \mbox{11 - Dense} & \mbox{$(B_\text{CL},64)$} & \mbox{Tanh} & \mbox{10}\\
    \mbox{12 - Dense} & \mbox{$(B_\text{CL},16)$} & \mbox{Tanh} & \mbox{11}\\
    \mbox{13 - Dense} & \mbox{$(B_\text{CL},N_y+1)$} & \mbox{$\text{Softmax}$} & \mbox{12}\\
        \bottomrule
          \end{tabular}
          }
       \subfloat[\label{tab:NN_CL_hyper}]{
       \scriptsize
  \begin{tabular}{ll}
    \toprule
            \mbox{Convolution kernel size:} & \mbox{$25,13,7$}\\
           \mbox{Dropout rate:} & \mbox{$5\%$}\\
\mbox{Weight initializer:} & \mbox{Xavier}\\
          \mbox{$L^2$ regularization rate:} & \mbox{$\lambda_\text{CL}=10^{-3}$}\\
\mbox{Optimizer:} & \mbox{Adam}\\
\mbox{Batch size:} & \mbox{$B_\text{CL}=32$}\\
\mbox{Initial learning rate:} & \mbox{$\eta_\text{CL}=\lbrace10^{-3},10^{-4}\rbrace$}\\
            \mbox{Allowed epochs:} & \mbox{$250$}\\
            \mbox{Learning schedule:} & \mbox{$\frac{4}{5}$ cosine decay}\\
            \mbox{Weight decay:} & \mbox{$0.05$}\\
\mbox{Early stop patience:} & \mbox{15 epochs}\\
\mbox{Train-val split:} & \mbox{$80:20$}\\
    \bottomrule
  \end{tabular}
  }
\end{table}

\begin{table}[h!]
\caption{$\text{N\hspace{-1px}N}_\text{RG}$ - (a) employed architecture, and (b) selected hyperparameters and training options.}
      \centering
       \subfloat[\label{tab:NN_RG_arch}]{
       \scriptsize
\begin{tabular}{llll}
    \toprule
    \mbox{Layer} & \mbox{Output shape} &  \mbox{Activation} & \mbox{Input}\\
    \toprule
    \mbox{0 - Input} & \mbox{$(B_\text{RG},L,N_u)$} & \mbox{None} & \mbox{--}\\
    \mbox{1 - Conv1D} & \mbox{$(B_\text{RG},L,32)$} & \mbox{Tanh} & \mbox{0}\\
    \mbox{2 - MaxPooling1D} & \mbox{$(B_\text{RG},L/2,32)$} & \mbox{None} & \mbox{1}\\
    \mbox{3 - Dropout} & \mbox{$(B_\text{RG},L/2,32)$} & \mbox{None} & \mbox{2}\\
    \mbox{4 - Conv1D} & \mbox{$(B_\text{RG},L/2,64)$} & \mbox{Tanh} & \mbox{3}\\
    \mbox{5 - MaxPooling1D} & \mbox{$(B_\text{RG},L/4,64)$} & \mbox{None} & \mbox{4}\\
    \mbox{6 - Dropout} & \mbox{$(B_\text{RG},L/4,64)$} & \mbox{None} & \mbox{5}\\
    \mbox{7 - Conv1D} & \mbox{$(B_\text{RG},L/4,32)$} & \mbox{Tanh} & \mbox{6}\\
    \mbox{8 - MaxPooling1D} & \mbox{$(B_\text{RG},L/8,32)$} & \mbox{None} & \mbox{7}\\
    \mbox{9 - Dropout} & \mbox{$(B_\text{RG},L/8,32)$} & \mbox{None} & \mbox{8}\\
    \mbox{10 - Flatten} & \mbox{$(B_\text{RG},4L)$} & \mbox{None} & \mbox{9}\\
    \mbox{11 - Dense} & \mbox{$(B_\text{RG},64)$} & \mbox{Tanh} & \mbox{10}\\
    \mbox{12 - Dense} & \mbox{$(B_\text{RG},16)$} & \mbox{Tanh} & \mbox{11}\\
    \mbox{13 - Dense} & \mbox{$(B_\text{RG},1)$} & \mbox{None} & \mbox{12}\\
        \bottomrule
        \end{tabular}
        }
        \subfloat[\label{tab:NN_RG_hyper}]{
       \scriptsize
  \begin{tabular}{ll}
    \toprule
          \mbox{Convolution kernel size:} & \mbox{$25,13,7$}\\
           \mbox{Dropout rate:} & \mbox{$10\%$}\\
\mbox{Weight initializer:} & \mbox{Xavier}\\
          \mbox{$L^2$ regularization rate:} & \mbox{$\lambda_\text{RG}=10^{-3}$}\\
\mbox{Optimizer:} & \mbox{Adam}\\
\mbox{Batch size:} & \mbox{$B_\text{RG}=32$}\\
\mbox{Initial learning rate:} & \mbox{$\eta_\text{RG}=\lbrace10^{-3},10^{-4}\rbrace$}\\
            \mbox{Allowed epochs:} & \mbox{$250$}\\
            \mbox{Learning schedule:} & \mbox{$\frac{4}{5}$ cosine decay}\\
            \mbox{Weight decay:} & \mbox{$0.05$}\\
\mbox{Early stop patience:} & \mbox{15 epochs}\\
\mbox{Train-val split:} & \mbox{$80:20$}\\
    \bottomrule
  \end{tabular}
  }
\end{table}

Using the Xavier's weight initialization~\cite{art:Glorot}, $\text{N\hspace{-1px}N}_\text{CL}$ and $\text{N\hspace{-1px}N}_\text{RG}$ are trained by minimizing the following loss functions, respectively:
\begin{align}
&\mathcal{L}^R_\text{CL}(\boldsymbol{\Theta}_\text{CL},\mathcal{D}_\text{CL}) = \mathcal{L}_\text{CL}(\boldsymbol{\Theta}_\text{CL},\mathcal{D}_\text{CL}) + \lambda_\text{CL}\lVert \boldsymbol{\Theta}_\text{CL} \rVert_{2}^2,\\
&\mathcal{L}^R_\text{RG}(\boldsymbol{\Theta}_\text{RG},\mathcal{D}_\text{RG}) = \mathcal{L}_\text{RG}(\boldsymbol{\Theta}_\text{RG},\mathcal{D}_\text{RG}) + \lambda_\text{RG}\lVert \boldsymbol{\Theta}_\text{RG} \rVert_{2}^2,
\label{eq:regularized_losses}
\end{align}
where $\lambda_\text{CL}$ and $\lambda_\text{RG}$ denote the $L^2$ regularization rate over the relative model parameters $\boldsymbol{\Theta}_\text{CL}$ and $\boldsymbol{\Theta}_\text{RG}$. The loss functions $\mathcal{L}^R_\text{CL}$ and $\mathcal{L}^R_\text{RG}$ are minimized using the first-order stochastic gradient descent optimizer Adam~\cite{art:Adam}, for a maximum of 250 allowed epochs. The corresponding learning rates $\eta_\text{CL}$ and $\eta_\text{RG}$ are initially set to $\lbrace10^{-3},10^{-4}\rbrace$, and decreased for $4/5$ of the allowed training steps using a cosine decay schedule with weight decay equal to $0.05$. The optimization is carried out considering an 80:20 splitting ratio of the dataset for training and validation purposes, with $20\%$ of the data randomly taken and set aside to monitor the learning process. We use an early stopping strategy to interrupt learning, whenever the loss function value attained on the validation set does not decrease for a prescribed number of patience epochs in a row. The hyperparameters and training options for $\text{N\hspace{-1px}N}_\text{CL}$ and for $\text{N\hspace{-1px}N}_\text{RG}$ are reported in \tab\ref{tab:NN_CL_hyper} and in \tab\ref{tab:NN_RG_hyper}, respectively. 

\end{document}